\documentclass[11pt,reqno,oneside]{amsart}

\usepackage{textcomp}
\usepackage{latexsym}
\usepackage{mathrsfs}
\usepackage{graphicx}
\usepackage[T1]{fontenc}
\usepackage[latin1]{inputenc}
\usepackage[english]{babel}
\usepackage[all]{xy} 
\usepackage{amsmath, amsfonts, amssymb, amscd, amsthm, enumerate}

\theoremstyle{plain}
\newtheorem{theo}{Theorem}[section]

\newtheorem{prop}[theo]{Proposizione}

\theoremstyle{definition}

\newenvironment{pf}{\noindent{\it Proof. }}{$\square$\par\medskip}

\theoremstyle{plain}
\newtheorem{lemma}[theo]{Lemma}
\newtheorem{theorem}[theo]{Theorem}

\newtheorem{proposition}[theo]{Proposition}

\theoremstyle{definition}

\newtheorem{remark}[theo]{Remark}
\newtheorem{defin}[theo]{Definition}
\newtheorem{example}[theo]{Example}

\newcommand{\beq}{\begin{equation}}
\newcommand{\eeq}{\end{equation}}

\renewcommand{\d}{\delta}
\newcommand{\e}{\epsilon}
\newcommand{\f}{\varphi}
\newcommand{\g}{\gamma}

\renewcommand{\l}{\lambda}

\newcommand{\q}{\vartheta}
\renewcommand{\r}{\rho}
\newcommand{\s}{\sigma}
\renewcommand{\t}{\tau}

\newcommand{\G}{\Gamma}

\renewcommand{\S}{\Sigma}
\newcommand{\X}{\Xi}

\newcommand{\bC}{\mathbb{C}}

\newcommand{\bN}{\mathbb{N}}

\newcommand{\bP}{\mathbb{P}}

\newcommand{\gD}{\mathfrak{D}}

\newcommand{\gU}{\mathfrak{U}}

\newcommand{\sF}{\mathscr{F}}

\newcommand{\cF}{\mathcal{F}}

\newcommand{\cI}{\mathcal{I}}

\newcommand{\cN}{\mathcal{N}}
\newcommand{\cO}{\mathcal{O}}

\newcommand{\cT}{\mathcal{T}}

\newcommand{\cV}{\mathcal{V}}

\newcommand{\cX}{\mathcal{X}}

\newcommand{\sD}{\mathscr{D}}

\renewcommand{\square}{\kern1pt\vbox
{\hrule height 0.6pt\hbox{\vrule width 0.6pt\hskip 3pt
\vbox{\vskip 6pt}\hskip 3pt\vrule width 0.6pt}\hrule height0.6pt}\kern1pt}

\newcommand{\n}{\nabla}

\newcommand{\ot}{\otimes}

\newcommand{\be}{\begin{equation}}
\newcommand{\ee}{\end{equation}}

\def\<#1,#2>{\langle\,#1,\,#2\,\rangle}
\newcommand{\arr}{\begin{array}{rlll}}
\newcommand{\ea}{\end{array}}
\newcommand{\bea}{\begin{eqnarray}}
\newcommand{\eea}{\end{eqnarray}}
\newcommand{\bean}{\begin{eqnarray*}}
\newcommand{\eean}{\end{eqnarray*}}
\newcommand{\ac}{\`}

\def\sideremark#1{\ifvmode\leavevmode\fi\vadjust{
\vbox to0pt{\hbox to 0pt{\hskip\hsize\hskip1em
\vbox{\hsize3cm\tiny\raggedright\pretolerance10000
\noindent #1\hfill}\hss}\vbox to8pt{\vfil}\vss}}}

\newcounter{ssig}
\newcounter{ttig}
\begin{document}

\title{Index theorems for couples of holomorphic self-maps}
\author{Paolo Arcangeli}
\address{Dipartimento di Matematica G. Castelnuovo, Universit\ac a ``Sapienza'' di Roma, Piazzale A. Moro 5, 00185, Roma, Italy}
\email{arcangeli@mat.uniroma1.it}
\begin{abstract}
Let $M$ be a $n$-dimensional complex manifold and $f,g:M\to M$ two distinct holomorphic self-maps. Suppose that $f$ and $g$ coincide on a globally irreducible compact hypersurface $S\subset M$. We show that if one of the two maps is a local biholomorphism around $S'=S-\text{Sing}(S)$ and, if needed, $S'$ sits into $M$ in a particular nice way, then it is possible to define a $1$-dimensional holomorphic (possibly singular) foliation on $S'$ and partial holomorphic connections on certain holomorphic vector bundles on $S'$. As a consequence, we are able to localize suitable characteristic classes and thus to get index theorems.\par
\bigskip
\noindent \textbf{\keywordsname :} Index theorem; Holomorphic foliation; Holomorphic connection;  Couple of holomorphic self-maps; Residue.
\end{abstract}
\maketitle

\section*{Introduction}
This paper is deeply inspired to the various works about index theorems for holomorphic self-maps and holomorphic foliations. Our goal here is to prove index theorems for couples $(f,g)$ of holomorphic self-maps coinciding on a positive dimensional set.\par
A first example of index theorem for holomorphic self-maps is the classical holomorphic Lefschetz fixed-point formula (see for example \cite[Ch.3,Sec.4]{GH1978}) which regards maps $f:M\to M$ having isolated fixed points, with $M$ a compact complex manifold. Anyway, we are mostly inspired by index theorems concerning self-maps having a positive dimensional fixed-points set, like for example the one in \cite{ab2001}. In his paper Abate obtained a complete generalization to two complex variables of the classical Leau-Fatou flower theorem for maps tangent to the identity and a key ingredient in the proof was an index theorem for holomorphic self-maps on complex surfaces $M$ fixing pointwise a non-singular compact complex curve $S\subset M$. This theorem was inspired by the Camacho-Sad index theorem for invariant leaves of possibly singular holomorphic foliations on complex surfaces (see \cite[Appendix]{CS1982}), later generalized to possibly singular leaves first (see \cite{su1995}) and then to arbitrary dimension of the ambient complex manifold, the foliation and the leaves (see \cite{le1991}, \cite{LS1995}, \cite{LS1999} or \cite{Su1998} for a complete treatment). Similarly, a first generalization of Abate index theorem was made in \cite{BT2002} (see also \cite{br2003}) where the authors assume $S$ to be possibly singular. However a large generalization to any dimension of $M$ and codimension of the possibly singular $S\subset M$ was made by Abate-Bracci-Tovena in \cite{ABT2004}, where they proved even other index theorems.\par
In this paper we do a step further because we replace the single holomorphic self-map $f:M\to M$ pointwise fixing an analytic sub-variety $S\subset M$ with a couple $(f,g)$ of distinct holomorphic self-maps coinciding on $S$ (which, for simplicity, we assume of codimension $1$). Clearly if one considers the couple $(f,\text{Id}_M)$ then falls back in the \cite{ABT2004} case. Briefly, we show that assuming some hypotheses on the couple $(f,g)$ and possibly on the hypersurface $S$ one can define a $1$-dimensional holomorphic foliation on $S'=S-\text{Sing}(S)$ and certain partial holomorphic connections (outside a `singular set') on suitable holomorphic vector bundles on $S'$. As a consequence one can use the Lehmann-Suwa machinery (see \cite{Su1998} or \cite{BSS2009} for a systematic exposition) in order to gain index theorems.\par
Our index theorems generalize the ones in \cite{ABT2004} and may be seen as versions for couples of holomorphic self-maps of some main index theorems of foliation theory. To be precise, we get new versions of the Baum-Bott index theorem (see \cite[Th.1.]{BB1972} or \cite[Th.III.7.6.]{Su1998}), of the above cited Camacho-Sad index theorem and of the Lehmann-Suwa (or variation) theorem (see \cite{LS1999} or \cite[Th.IV.5.6.]{Su1998}). We point out that for the last two cited index theorems of foliation theory one needs to have a foliation defined around $S$ leaving it invariant, while we are able to define foliations on $S$ only. However, the foliations we define extend naturally to a suitable infinitesimal neighborhood of the sub-variety, and this allows to localize certain characteristic classes producing our index theorems (see also \cite{ABT2008} for a general explanation).\par
\smallskip
The plan of the paper and our main theorems are the following. In Section \ref{ord_can_sec} we introduce the `order of coincidence of $(f,g)$ along $S$', which is a positive integer denoted by $\nu_{f,g}$, and then we define the `canonical section associated to $(f,g)$'
$$
\sD_{f,g}:N_{S'}^{\ot\nu_{f,g}}\longrightarrow TM|_{S'}
$$
on the regular part $S'$ of $S$. In order to define the canonical section we need to assume that one of the two maps is a local biholomorphism on a whole neighborhood of $S'$ and we will make this assumption for the rest of the paper. In the following Section \ref{foliations} we introduce two hypotheses under which the canonical section associated to $(f,g)$ induces a $1$-dimensional holomorphic foliation on $S'$ (or more than one in some cases), denoted here by
$$
\sD:N_{S'}^{\ot\nu_{f,g}}\longrightarrow TS'.
$$
One hypothesis regards the couple $(f,g)$, which is said `tangential along $S$' when it is satisfied. The other concerns the way $S'$ sits into $M$ and in this case we say that `$S'$ splits into $M$'. Subsequently, we prove in Section \ref{BaumBott} that when $S$ is non-singular and one of the two cited hypotheses is verified one can define in a canonical way a partial holomorphic connection (in the sense of Bott \cite{bo1967}) $\d_{\sD}^{bb}$ on the normal bundle $N_{\sD}$ of the foliation. As a consequence we get an index theorem which we state here in a simplified version.
\begin{theorem}\label{th1}
Let $M$ be a $n$-dimensional complex manifold, $S\subset M$ a non-singular compact connected complex hypersurface in $M$ and $(f,g)$ a couple of holomorphic self-maps on $M$  such that $f|_S=g|_S$ and $g$ is a local biholomorphism on a neighborhood of $S$. Assume that
\begin{itemize}
\item[(i)] $(f,g)$ is tangential along $S$ 
\end{itemize}
or that
\begin{itemize}
\item[(ii)] $S$ splits into $M$
\end{itemize}
and let $\sD$ be the foliation on $S$ induced by the canonical section $\sD_{f,g}$. Assume $\sD\neq 0$ and let $\text{Sing}(\sD)=\sqcup_{\l}\S_{\l}$ be the decomposition in connected components of the singular set of the foliation. \par
Then for any symmetric homogeneous polynomial $\f\in\bC[z_1,\dots,z_{n-2}]$ of degree $n-1$ there exist complex numbers $\text{Res}_{\f}(\sD;TS-N_S^{\ot\nu};\S_{\l})$ such that
$$
\sum_{\l}\text{Res}_{\f}\left(\sD;TS-N_S^{\ot\nu};\S_{\l}\right)=\int_S\f\left(TS-N_S^{\ot\nu}\right).
$$
\end{theorem}   
We conclude by computing the residues appearing in Theorem \ref{th1} at isolated singular points of the foliation. In Section \ref{loc_ext} we show that if $(f,g)$ is tangential along $S$ then one can locally extend the foliation $\sD$ on $S'$ in a `canonical way' and control in a sense the differences among canonical local extensions. We see moreover that if $(f,g)$ is not tangential along $S$ one can still do this but assuming a hypothesis on $S'$ stronger than the splitting property. When this condition is satisfied we say that `$S'$ is comfortably embedded into $M$'. Finally, in the last two sections we reap the benefits of Section \ref{loc_ext}. Indeed in Section \ref{CamachoSad} we show that these canonical local extensions of $\sD$ are good enough to define (outside a `singular set') a partial holomorphic connection $\d_{\sD}^{cs}$ on the normal bundle $N_{S'}$ of $S'$ in $M$. The result is the following theorem.
\begin{theorem}\label{th2}
Let $M$ be a $n$-dimensional complex manifold, $S\subset M$ a globally irreducible compact complex hypersurface and $(f,g)$ a couple of holomorphic self-maps on $M$  such that $f|_S=g|_S$ and $g$ is a local biholomorphism on a neighborhood of $S'$. Assume that
\begin{itemize}
\item[(i)] $(f,g)$ is tangential along $S$ 
\end{itemize}
or that
\begin{itemize}
\item[(ii)] $S'$ is comfortably embedded into $M$
\end{itemize}
and let $\sD$ be the foliation on $S'$ induced by the canonical section $\sD_{f,g}$. Assume $\sD\neq 0$ and let $\text{Sing}(S)\cup\text{Sing}(\sD)=\sqcup_{\l}\S_{\l}$ be the decomposition in connected components of the singular set $\text{Sing}(S)\cup\text{Sing}(\sD)$. \par
Then there exist complex numbers $\text{Res}(\sD;S;\S_{\l})$ such that
$$
\sum_{\l}\text{Res}\left(\sD;S;\S_{\l}\right)=\int_S c_1^{n-1}\left([S]\right),
$$
where $c_1([S])$ denotes the first Chern class of the line bundle $[S]$ on $M$ canonically associated to $S$.
\end{theorem}
Again, we conclude by computing the residues appearing in Theorem \ref{th2} at isolated singular points. Similarly, in Section \ref{LehmannSuwa} we show that if $(f,g)$ is tangential along $S$ and $\nu_{f,g}>1$  then the canonical local extensions of $\sD$ are good enough to define a partial holomorphic connection $\d_{\sD}^{ls}$ on the normal bundle $N_{\sD}^M$ of the foliation respect to the ambient tangent bundle $TM$ (restricted to $S'$). It follows the last index theorem.
\begin{theorem}\label{th3}
Let $M$ be a $n$-dimensional complex manifold, $S\subset M$ a globally irreducible compact complex hypersurface and $(f,g)$ a couple of holomorphic self-maps on $M$  such that $f|_S=g|_S$ and $g$ is a local biholomorphism on a neighborhood of $S'$. Suppose $(f,g)$ tangential along $S$ and $\nu_{f,g}>1$ and let $\sD$ be the foliation on $S'$ induced by the canonical section. Assume $\sD\neq 0$ and let $\text{Sing}(S)\cup\text{Sing}(\sD)=\sqcup_{\l}\S_{\l}$ be the decomposition in connected components of the singular set $\text{Sing}(S)\cup\text{Sing}(\sD)$. \par
Then for any symmetric homogeneous polynomial $\f\in\bC[z_1,\dots,z_{n-1}]$ of degree $n-1$ there exist complex numbers $\text{Res}_{\f}(\sD;TM|_S-[S]|_S^{\ot\nu};\S_{\l})$ such that
$$
\sum_{\l}\text{Res}\left(\sD;\left.TM\right|_S-\left.[S]\right|_S^{\ot\nu};\S_{\l}\right)=\int_S \f\left(TM-[S]^{\ot\nu}\right).
$$
\end{theorem}
As for the other index theorems we end by computing the residues appearing in Theorem \ref{th3} at isolated singular points.\par
\bigskip
\section{The order of coincidence and the canonical section}\label{ord_can_sec}\setcounter{equation}{0}
Let $M$ be a $n$-dimensional complex manifold and $S\subset M$ a (possibly singular) globally irreducible complex hypersurface. From now on we will denote by $\cO_M$ the sheaf of germs of holomorphic functions on $M$ and by $\cI_S$ the sub-sheaf of germs of functions vanishing on $S$, which is a coherent $\cO_M$-module. For the sake of simplicity we shall use the same symbol to denote both a germ at some point and any representative defined in a neighborhood of the point. Recall that the sheaf of germs of holomorphic functions on $S$ is by definition $\cO_S=\cO_M/\cI_S$. We will denote by $TM$ the holomorphic tangent bundle of $M$ and, in case $S$ is non-singular, by $TS$ the holomorphic tangent bundle of $S$ and by $N_S$ its holomorphic normal bundle in $M$, which is a line bundle. The corresponding sheaves of germs of holomorphic sections will be denoted respectively by $\cT_M$, $\cT_S$ and $\cN_S$. Lastly, we will denote by $\text{End}^2_S(M)$ the set of couples $(f,g)$ of  distinct (germs about $S$ of) holomorphic self-maps of $M$ {\it coinciding on $S$}. In other words, if $(f,g)\in\text{End}^2_S(M)$ then $f,g:M\to M$ are holomorphic, $f\neq g$, and $f|_S\equiv g|_S$.\par
The goal of this first section is to introduce the concept of `order of coincidence' of a couple $(f,g)$ and to define the `canonical section' associated to a couple with a certain property, which is a holomorphic section of a suitable holomorphic vector bundle on $S'=S-\text{Sing}(S)$.\par
\smallskip
So let $(f,g)\in\text{End}^2_S(M)$ be a given couple and fix a point $p\in S$. Observe that for every germ $h\in\cO_{M,f(p)}$ we have the well-defined germ $h\circ f-h\circ g\in\cI_{S,p}\subset\cO_{M,p}$, so we can give the following definition.\par
\begin{defin}\label{def_order}
For any $h\in\cO_{M,f(p)}$, we define $\nu_{f,g}^p(h)$ to be the constant
$$
\nu_{f,g}^p(h)=\text{max}\{\nu\in\bN\text{ s.t. }h\circ f-h\circ g\in\cI^{\nu}_{S,p}\}.
$$
The {\it order of coincidence of $(f,g)$ at $p$ along $S$} is then 
$$
\nu_{f,g}^p=\text{min}\{\nu_{f,g}^p(h),\text{ for }h\in\cO_{M,f(p)}\}.
$$
\end{defin}
If $(w^1,\dots,w^n)$ is any local coordinates system at $f(p)=g(p)$ then for every $h\in\cO_{M,f(p)}$ we have the fundamental relations of germs at $p$
\begin{align}\label{fundamentals}
h\circ f-h\circ g & =\sum_{j=1}^n\left(f^j-g^j\right)\frac{\partial h}{\partial w^j}\circ g \left(\text{mod }\cI^{2\nu_{f,g}^p}_{S,p}\right)= \notag \\
 & =\sum_{j=1}^n\left(f^j-g^j\right)\frac{\partial h}{\partial w^j}\circ f \left(\text{mod }\cI^{2\nu_{f,g}^p}_{S,p}\right)
\end{align}
where $f^j=w^j\circ f$ and $g^j=w^j\circ g$. In fact, by Definition \ref{def_order}  
\begin{align*}
f^j-g^j=s^j& \in\cI^{\nu_{f,g}^p}_{S,p},\qquad\qquad j=1,\dots,n\\
\frac{\partial h}{\partial w^j}\circ f-\frac{\partial h}{\partial w^j}\circ g & \in\cI^{\nu_{f,g}^p}_{S,p},\qquad\qquad j=1,\dots,n,
\end{align*}
then
\begin{align*}
h\circ f & -h\circ g=h(g^1+s^1,\dots,g^n+s^n)-h(g^1,\dots,g^n)=\\
& =\sum_{j=1}^n s^j\left.\frac{\partial h}{\partial w^j}\right|_{(g^1,\dots,g^n)}+\frac{1}{2}\sum_{j,k=1}^n s^js^k\left.\frac{\partial^2 h}{\partial w^j\partial w^k}\right|_{(g^1,\dots,g^n)}+\cdots=\\
& =\sum_{j=1}^n(f^j-g^j)\frac{\partial h}{\partial w^j}\circ g\left(\text{mod }\cI^{2\nu_{f,g}^p}_{S,p}\right)=\\
& =\sum_{j=1}^n(f^j-g^j)\frac{\partial h}{\partial w^j}\circ f\left(\text{mod }\cI^{2\nu_{f,g}^p}_{S,p}\right),
\end{align*}
where the second row is given by the Taylor expansion of $h$ at $f(p)=g(p)$. \par
\begin{lemma}\label{lemma1}
Let $p\in S$. Then
\begin{enumerate}
\item[i)]if $(w^1,\dots,w^n)$ is any set of local coordinates at $f(p)$ then
$$
\nu_{f,g}^p=\text{min}\left\{\nu_{f,g}^p(w^1),\dots,\nu_{f,g}^p(w^n)\right\}.
$$
\item[ii)]for any $h\in\cO_{M,f(p)}$ the function
$$
q\to \nu_{f,g}^q(h)
$$
is constant in a neighborhood of $p$ in $S$. 
\item[iii)]the function
$$
p\to \nu_{f,g}^p
$$
is constant on $S$.
\end{enumerate}
\end{lemma}
\begin{pf}
i) By Definition \ref{def_order} it follows that $\nu_{f,g}^p$$\leq$$\text{min}\left\{\nu_{f,g}^p(w^1),\dots,\nu_{f,g}^p(w^n)\right\}$ and by (\ref{fundamentals}) it follows the reverse inequality.\par
ii) Let $h\in\cO_{M,f(p)}$ and $\{\g^1,\dots,\g^t\}$ be a set of generators of $\cI_{S,p}$, then
$$
h\circ f-h\circ g=\sum_{|I|=\nu_{f,g}^p(h)}\g^Ic_I,
$$
where $I=(i_1,\dots,i_t)\in\bN^t$ is a multi-index, $|I|=i_1+\cdots+i_t$, $\g^I=(\g^1)^{i_1}\cdots (\g^t)^{i_t}$ and $c_I\in\cO_{M,p}$. This equality clearly holds in a neighborhood of $p$ by definition of germ. Moreover, since $\cI_S$ is a coherent sheaf the $\g$'s are generators of $\cI_{S,q}$ for points $q$ near enough to $p$. Finally, observe that by definition of $\nu_{f,g}^p(h)$ there is an index $I_0$ such that $c_{I_0}\notin \cI_{S,p}$, then $c_{I_0}\notin \cI_{S,q}$ for all $q\in S$ close enough to $p$. All these facts make the assertion follows easily.\par
iii) By i) and ii) the function $p\to\nu_{f,g}^p$ is locally constant on $S$. Since $S$ is connected then it is constant.
\end{pf}
As a consequence of $iii)$ the following definition makes sense.
\begin{defin}
The {\it order of coincidence of $(f,g)$ along $S$} is the constant
$$
\nu_{f,g}=\nu_{f,g}^p
$$
for any point $p\in S$.
\end{defin}
\par\smallskip
Now assume $S$ non-singular in the following and let $(f,g)\in\text{End}^2_S(M)$ be a couple in which one of the two maps, say $g$, is a local biholomorphism if restricted to an open neighborhood of $S$ in $M$. This means that there exists an open neighborhood $W\subset M$ containing $S$ such that $g|_W:W\to M$ is a local biholomorphism, or equivalently that $\text{d}g|_p:T_pM\longrightarrow T_{g(p)}M $ is an isomorphism of complex vector spaces for every $p\in S$. As $f|_S\equiv g|_S$ there is a well-defined morphism of holomorphic vector bundles
\begin{equation}\label{pre_sec}
\text{d}f-\text{d}g: N_S \longrightarrow \left.(f^* TM)\right|_S,
\end{equation}
which on the fibers is given by $[v] \to [\text{d}f|_p(v)-\text{d}g|_p(v)]$, for any $[v]\in T_pM$ and $p\in S$. Observe that while in general $f^* TM\neq g^* TM$ obviously by the assumptions $(f^* TM)|_S=(g^* TM)|_S$. By the property of $g$ we also have the isomorphism of holomorphic vector bundles 
$$
\left.\text{d}g\right|_S:TM|_S\longrightarrow \left(f^* TM\right)|_S
$$
thus composing its inverse with (\ref{pre_sec}) we obtain the morphism
\begin{equation}\label{sec}
\left.\text{d}g\right|_S^{-1}\circ (\text{d}f-\text{d}g):N_S\longrightarrow TM|_S,
\end{equation}
which can also be seen as a holomorphic section of $N_S^*\ot TM|_S$. We want to express (\ref{sec}) in a local frame but first we introduce some general terminology about local charts on $M$.
\begin{defin}\label{def_coords}
Let $M$ be a $n$-dimensional complex manifold and $S\subset M$ a complex sub-manifold of any codimension $k$ ($0<k<n$). We say that a local holomorphic chart $(U,z)=(U,z^1,\dots,z^n)$ is {\it adapted to $S$} if $U\cap S=\emptyset$ or if $U\cap S=\{z^1=\cdots=z^k=0\}$. Equivalently, we can say that the local coordinates are adapted to $S$. If moreover the coordinates are centered at $p\in U\cap S$ we say that the chart is (or the coordinates are) {\it adapted to $S$ at $p$}. We call an atlas $\gU$ of $M$ of local adapted charts an {\it atlas adapted to $S$}.\par
Let $W\subset M$ be an open neighborhood of $S$ and $\f:W\to M$ a local biholomorphism. If $(U,z)$ is a local chart such that $\left. \f\right|_U$ is a biholomorphism onto its image we can consider the coordinates $w=(w^1=z^1\circ \f^{-1},\dots,w^n=z^n\circ \f^{-1})$ on $g(U)$. If $(U,z)$ is also adapted to $S$ (at $p$)  we say that the coordinates $z$ are (or the chart $(U,z)$ is) {\it adapted to $(\f,S)$ (at $p$)} and the coordinates $w=z\circ \f^{-1}$ are called {\it special}. We call an atlas $\gU$ of $M$ of local charts adapted to $(\f,S)$ an {\it atlas adapted to $(\f,S)$}.
\end{defin}
Observe that if $S\subset M$ is a hypersurface and $(U,z)$ is a local chart adapted to it such that $U\cap S\neq\emptyset$ then
$$
\left\{\left.\frac{\partial}{\partial z^2}\right|_S,\dots,\left.\frac{\partial}{\partial z^n}\right|_S\right\}
$$
is a local holomorphic frame for $TS$, while if $\pi:TM|_S\to N_S$ is the obvious projection then
$$
\partial z^1=\pi\left(\left.\frac{\partial}{\partial z^1}\right|_S\right)
$$
is a local holomorphic generator for $N_S$ and we denote with $\partial^* z^1$ its dual (which is a local holomorphic generator for $N_S^*$).\par
In order to write (\ref{sec}) locally we calculate (\ref{pre_sec}) first. So let $(U,z)$ be a local chart adapted to $S$ at a point $p\in S$ and $(V,w)$ be {\it any} local chart at $f(p)$. Then
$$
\left\{\partial^*z^1\ot f^*\left.\frac{\partial}{\partial w^1}\right|_S,\dots,\partial^*z^1\ot f^*\left.\frac{\partial}{\partial w^n}\right|_S\right\}
$$
is a local holomorphic frame for $N_S^*\ot (f^*TM)|_S$ on the neighborhood $U\cap f^{-1}(V)\cap S$ of $p$, while
$$
\left\{\partial^*z^1\ot \left.\frac{\partial}{\partial z^1}\right|_S,\dots,\partial^*z^1\ot \left. \frac{\partial}{\partial z^n}\right|_S\right\}
$$
is a local holomorphic frame for $N_S^*\ot TM|_S$ on $U\cap S$. Observe that even if in general $f^*\frac{\partial}{\partial w^j}\neq g^*\frac{\partial}{\partial w^j}$, their restrictions to $S$ are equal by the assumptions. If we denote $f^j=w^j\circ f$ and $g^j=w^j\circ g$ then clearly $f^j-g^j\in\cI_{S,p}$ for every $j$ by the hypothesis. Since the coordinates $z$ are adapted to $S$ there exist germs $h^j\in\cO_{M,p}$ such that $f^j-g^j=h^j z^1$, for $j=1,\dots,n$. A trivial calculation shows that (\ref{pre_sec}) is locally given by
\begin{equation}\label{loc_pre_sec}
\sum_{j=1}^n\left.\frac{\partial(f^j-g^j)}{\partial z^1}\right|_S \partial^*z^1\ot \left.f^*\frac{\partial}{\partial w^j}\right|_S=\sum_{j=1}^n\left.h^j\right|_S \partial^*z^1\ot \left.f^*\frac{\partial}{\partial w^j}\right|_S.
\end{equation}
If $(U,z)$ is adapted to $(g,S)$ and we take the associated special coordinates $w$ on $V=g(U)$ we can easily compute (\ref{sec}). In fact, with this choice of coordinates 
$$
\left. f^*\frac{\partial}{\partial w^j}\right|_q=\left. g^*\frac{\partial}{\partial w^j}\right|_q=\left.\frac{\partial}{\partial w^j}\right|_{g(q)}=\left.\text{d}g\right|_q\left(\left.\frac{\partial}{\partial z^j}\right|_q \right)
$$
for every $q\in U\cap S$ and $j=1,\dots,n$, so by (\ref{loc_pre_sec}) it follows that (\ref{sec}) is locally 
\begin{equation}\label{loc_sec}
\sum_{j=1}^n\left.h^j\right|_S \partial^*z^1\ot \left.\frac{\partial}{\partial z^j}\right|_S.
\end{equation}
By Lemma \ref{lemma1} and (\ref{loc_sec}) it follows that the morphism (\ref{sec}) vanishes identically on $S$ if and only if $\nu_{f,g} >1$, then it would not be a good `canonical section'. The idea is to take local higher order derivatives (respect to $z^1$) of ``$f-g$'', not vanishing on $S$. To be more precise let $p\in S$ be any point, $(U,z)$ a local chart adapted to $(g,S)$ at $p$ and take on $g(U)$ the special coordinates $(w^1,\dots,w^n)$. Since the coordinates $z$ are adapted to $S$ there exist suitable germs $h^j\in\cO_{M,p}$ such that
\begin{equation}\label{germs_h}
f^j-g^j=w^j\circ f-w^j\circ g=h^j(z^1)^{\nu_{f,g}}
\end{equation}
for $j=1,\dots,n$. Let define
\begin{equation}\label{XXX}
\gD_{f,g}:=\sum_{j=1}^n h^j(\text{d}z^1)^{\nu_{f,g}}\ot\frac{\partial}{\partial z^j},
\end{equation}
where $(\text{d}z^1)^{\nu_{f,g}}=(\text{d}z^1)^{\ot\nu_{f,g}}$, which is a local holomorphic section of $(TM^{\ot\nu_{f,g}})^*\ot TM$ on $U$. Observe that it does not vanish identically on $S$ by Lemma \ref{lemma1}.\par
\begin{remark}
A priori the germs $h^j$ do not have representatives defined on $U$ but we can assume it (possibly shrinking $U$). From now on we will assume the $h^j$ to be defined on the whole $U$.
\end{remark}
We have the following fundamental proposition.
\begin{proposition}\label{gluing}
Let $(U,z)$, $(\hat{U},\hat{z})$ be two local chart adapted to $(g,S)$ at $p$ and $\gD_{f,g}$, $\hat{\gD}_{f,g}$ the corresponding local section of $(TM^{\ot\nu_{f,g}})^*\ot TM$ defined about $p$ as in (\ref{XXX}). Then 
$$
\hat{\gD}_{f,g}=\gD_{f,g}\quad\left(\text{mod }\cI_S\right)
$$
where they overlap.
\end{proposition}
\begin{pf}
In the following set $\nu=\nu_{f,g}$ for ease the notation. Since $z^1$ and $\hat{z}^1$ are both generators of $\cI_{S,p}$ it follows that
\begin{align}\label{alpha}
\hat{z}^1=a z^1&  &\text{and then}& & (\text{d}\hat{z}^1)^{\nu}=a^{\nu}(\text{d}z^1)^{\nu}\quad\left(\text{mod }\cI_S\right)
\end{align}
for some germ $a\in\cO_{M,p}^*$. Using (\ref{fundamentals}), (\ref{germs_h}) and (\ref{alpha}) we have that
\begin{align*}
\hat{h}^j a^{\nu}(z^1)^{\nu}=\sum_{k=1}^n h^k(z^1)^{\nu}\left(\frac{\partial\hat{w}^j}{\partial w^k}\circ g\right)\phantom{a}\left(\text{mod }\cI_S^{2\nu}\right),
\end{align*}
hence
\begin{equation}\label{htilde_j}
\hat{h}^j a^{\nu}=\sum_{k=1}^n h^k \frac{\partial \hat{z}^j}{\partial z^k}\phantom{a}\left(\text{mod }\cI_S^{\nu}\right),\qquad\quad j=1,\dots,n,
\end{equation}
since one can easily check that $\frac{\partial\hat{w}^j}{\partial w^k}\circ g=\frac{\partial \hat{z}^j}{\partial z^k}$. In particular for $j=1$ 
\begin{equation}\label{htilde_1}
\hat{h}^1 a^{\nu}=h^1 a \phantom{a}\left(\text{mod }\cI_S\right).
\end{equation}
An easy computation shows that
\begin{equation}\label{deriv1}
\frac{\partial}{\partial\hat{z}^1}=\frac{1}{a}\frac{\partial}{\partial z^1}+\sum_{k=2}^n\frac{\partial z^k}{\partial\hat{z}^1}\frac{\partial}{\partial z^k}\phantom{a}\left(\text{mod }\cI_S\right),
\end{equation}
\begin{equation}\label{deriv2}
\frac{\partial}{\partial\hat{z}^j}=\sum_{k=2}^n\frac{\partial z^k}{\partial\hat{z}^j}\frac{\partial}{\partial z^k}\phantom{a}\left(\text{mod }\cI_S\right),\qquad j=2,\dots,n.
\end{equation}
Then using (\ref{alpha}), (\ref{htilde_j}), (\ref{htilde_1}), (\ref{deriv1}) and (\ref{deriv2}) it follows that
\begin{align}\label{pre_end}
\hat{\gD}_{f,g} & =\sum_{j=1}^n\hat{h}^j(\text{d}\hat{z}^1)^{\nu}\ot\frac{\partial}{\partial\hat{z}^j}=\sum_{j=1}^n\hat{h}^j a^{\nu}(\text{d}z^1)^{\nu}\ot\frac{\partial}{\partial\hat{z}^j}\phantom{a}\left(\text{mod }\cI_S\right)= \notag \\
 & =a h^1 (\text{d}z^1)^{\nu}\ot\frac{\partial}{\partial\hat{z}^1}+\sum_{j=2}^n \sum_{r=1}^n h^r\frac{\partial \hat{z}^j}{\partial z^r}(\text{d}z^1)^{\nu}\ot\frac{\partial}{\partial\hat{z}^j}  \phantom{a}\left(\text{mod }\cI_S\right)= \notag \\
 & =h^1 (\text{d}z^1)^{\nu}\ot\frac{\partial}{\partial z^1}+a h^1\sum_{k=2}^n\frac{\partial z^k}{\partial\hat{z}^1}(\text{d}z^1)^{\nu}\ot \frac{\partial}{\partial z^k} +\notag \\
 & +\sum_{j=2}^n \sum_{r=1}^n \sum_{k=2}^n h^r\frac{\partial z^k}{\partial\hat{z}^j}\frac{\partial \hat{z}^j}{\partial z^r} (\text{d}z^1)^{\nu}\ot \frac{\partial}{\partial z^k}\phantom{a}\left(\text{mod }\cI_S\right) = \notag \\     
 & =h^1 (\text{d}z^1)^{\nu}\ot\frac{\partial}{\partial z^1}+h^1\sum_{k=2}^n\left[a\frac{\partial z^k}{\partial\hat{z}^1}+\sum_{j=2}^n\frac{\partial z^k}{\partial\hat{z}^j}\frac{\partial\hat{z}^j}{\partial z^1}\right](\text{d}z^1)^{\nu}\ot\frac{\partial}{\partial z^k}+ \notag \\
 & +\sum_{j=2}^n \sum_{r=2}^n \sum_{k=2}^n h^r\frac{\partial z^k}{\partial\hat{z}^j}\frac{\partial \hat{z}^j}{\partial z^r} (\text{d}z^1)^{\nu}\ot \frac{\partial}{\partial z^k}\phantom{a}\left(\text{mod }\cI_S\right) 
\end{align}
To conclude observe that for $k=2,\dots,n$
$$
0=\frac{\partial z^k}{\partial z^1}=\sum_{j=1}^n\frac{\partial z^k}{\partial\hat{z}^j}\frac{\partial\hat{z}^j}{\partial z^1}=a\frac{\partial z^k}{\partial \hat{z}^1}+\sum_{j=2}^n\frac{\partial z^k}{\partial\hat{z}^j}\frac{\partial\hat{z}^j}{\partial z^1}\phantom{a}\left(\text{mod }\cI_S\right)
$$
and that
$$
\d_{kr}=\frac{\partial z^k}{\partial z^r}=\sum_{j=1}^n\frac{\partial z^k}{\partial\hat{z}^j}\frac{\partial\hat{z}^j}{\partial z^r}=\sum_{j=2}^n\frac{\partial z^k}{\partial\hat{z}^j}\frac{\partial\hat{z}^j}{\partial z^r}\phantom{a}\left(\text{mod }\cI_S\right),
$$
for $r=2,\dots,n$. Putting these relations into (\ref{pre_end}) we have done.
\end{pf}
Thanks to Proposition \ref{gluing} we are now able to define the `canonical section'.
\begin{defin}
Let $M$ be a $n$-dimensional complex manifold, $S\subset M$ a non-singular connected complex hypersurface and $(f,g)\in\text{End}^2_S(M)$ a couple with $g$ a local biholomorphism if restricted to an open neighborhood of $S$.\par
The {\it canonical section associated to $(f,g)$} is the global holomorphic section $\sD_{f,g}$ of the holomorphic vector bundle $(N_S^{\ot\nu_{f,g}})^*\ot TM|_S$ obtained by gluing together {\it on $S$} the local $\gD_{f,g}$ defined in (\ref{XXX}). We can also think to it as a holomorphic section of $\text{Hom}(N_S^{\ot\nu_{f,g}},TM|_S)$, that is as a morphism $\sD_{f,g}:N_S^{\ot\nu_{f,g}}\to TM|_S$ of holomorphic vector bundles.
\end{defin}
Summing up, if $(z^1,\dots,z^n)$ are local coordinates adapted to $(g,S)$ then $\sD_{f,g}$ is locally defined by
$$
\sD_{f,g}\stackrel{loc.}{=}\sum_{j=1}^n\left.h^j\right|_S(\partial^* z^1)^{\nu}\ot\left.\frac{\partial}{\partial z^j}\right|_S,
$$
where the $h^j$ are the ones appearing in (\ref{germs_h}). Observe that $\sD_{f,g}$ is not identically vanishing on $S$ by construction.
\begin{defin}
A point $p\in S$ is a {\it singularity of  $\sD_{f,g}$} if the associated morphism $N_S^{\ot\nu_{f,g}}\to TM|_S$ is not injective in $p$, that is if  $\sD_{f,g}(p)=0$. We denote the {\it singular set of $\sD_{f,g}$} by $\text{Sing}(f,g)$.
\end{defin}
\begin{remark}
If we consider the couple $(f,\text{Id}_M)$ we are in the setting of \cite{ABT2004}. In this case every local adapted chart is clearly also $\text{Id}_{M}$-adapted, moreover Proposition \ref{gluing} turns out to be \cite[Prop.3.1.]{ABT2004} and the canonical section $\sD_{f,\text{Id}_M}$ is exactly the canonical section $X_f$ of \cite[Def.3.2.]{ABT2004}.
\end{remark}
\bigskip
\section{The canonical foliations}\label{foliations}\setcounter{equation}{0}
Let $M$ and $S$ be as in Section \ref{ord_can_sec}, with $S'=S-\text{Sing}(S)$ the regular part of $S$, and let $(f,g)\in\text{End}_S^2(M)$ be a couple whose order of coincidence is $\nu=\nu_{f,g}$ and in which $g$ is a local biholomorphism around $S'$ . As just seen, $(f,g)$ induces a canonical section $\sD_{f,g}:N_{S'}^{\ot\nu}\to TM|_{S'}$ but we would like to have a $1$-dimensional holomorphic (possibly singular) foliation on $S'$, that is an injective morphism 
$
\sF:\cF\to\cT_{S'}
$
of $\cO_{S'}$-modules where $\cF$ is a rank $1$ locally free $\cO_{S'}$-module. Recall that its (possibly) singular set is
$$
\text{Sing}(\sF)=\left\{x\in S'\text{ s.t. }(\cT_{S'}/\cF)_x \text{ is not a free $\cO_{S',x}$-module}\right\}
$$
Equivalently, a foliation on $S'$ can be described as a morphism 
$
\sF:F\to TS'
$
of holomorphic vector bundles on $S'$ where $F$ is a line bundle. In this case its (possibly) singular set could be described as
$$
\text{Sing}(\sF)=\left\{x\in S'\text{ s.t. } \sF(x)=0\right\}.
$$ 
\begin{defin}
The $\cO_{S'}$-module $\cF$ is called {\it tangent sheaf of $\sF$}, while the holomorphic vector bundle $F$ is the {\it tangent bundle of $\sF$}. Clearly $\cF$ is the sheaf of germs of holomorphic sections of $F$.
\end{defin}
In the current section we will discuss two conditions under which we have a foliation (or more foliations) on $S'$ induced by the canonical section $\sD_{f,g}$. The first is a condition on the couple $(f,g)$ while the second is on the way $S'$ is embedded into $M$.\par
\smallskip
Let $p\in S$ be any point. We have two induced ring homomorphisms given by pull-back of germs by $f$ or $g$,
$$
\cO_{M,f(p)}\stackrel{f^*_p,g^*_p}{\longrightarrow}\cO_{M,p}.   
$$
Clearly in general $f^*_p\neq g^*_p$ but since $f|_S =g|_S$ then $(f^*_p)^{-1}(\cI_{S,p})=(g^*_p)^{-1}(\cI_{S,p})$ and we denote this ideal of $\cO_{M,f(p)}$ by $I_{f(p)}$.\par
\begin{remark}\label{remark1}
Observe that $I_{f(p)}\supseteq \cI_{f(S),f(p)}$ but they are not the same in general. Anyway $I_{f(p)}$ may be seen as the stalk of a coherent module defined on an open neighborhood of $f(p)$, for every $p\in S$. In fact, by the hypothesis for any $p\in S$ there is an open neighborhood $U \subset M$ such that $g:U\to g(U)$ is a biholomorphism and then $f(U\cap S)=g(U\cap S)$ is an analytic sub-variety of $g(U)\subset M$. Hence $\cI_{f(U\cap S)}$ is a coherent $\cO_{g(U)}$-module and $I_{f(p)}=\cI_{f(U\cap S),f(p)}$. 
\end{remark}
\begin{defin}\label{def_tang}
Let $p\in S$. If 
$$
\text{min}\{\nu_{f,g}^p(h),\text{ for }h\in I_{f(p)}\}>\nu_{f,g}^p
$$
we say that $(f,g)$ is {\it tangential at $p\in S$}.\par
\end{defin}
\begin{lemma}\label{lemma_tang}
The following two statements are true:
\begin{enumerate}
\item[i)] Let $p\in S$ be any point and $\{\r^1,\dots,\r^k\}$ any set of generators of $I_{f(p)}$. Then $(f,g)$ is tangential at $p$ if and only if  
$$
\text{min}\{\nu_{f,g}^p(\r^1),\dots,\nu_{f,g}^p(\r^k)\}>\nu_{f,g}^p.
$$
\item[ii)] If $(f,g)$ is tangential at a point $p\in S$, then it is tangential at all the points of $S$.
\end{enumerate}
\end{lemma}
\begin{pf}
i) Let $h\in I_{f(p)}$, then $h=h_1\r^1+\dots h_k\r^k$ for some $h_1,\dots,h_k\in\cO_{M,f(p)}$. We can easily check that
$$
h\circ f-h\circ g=\sum_{j=1}^k (h_j\circ f)\left(\r^j\circ f-\r^j\circ g\right)+\sum_{j=1}^k (\r^j\circ g)\left(h_j\circ f-h_j\circ g\right).
$$
The $j$-th term of the first sum is in $\cI_{S,p}^{\nu_{f,g}^p(\r^j)}$ and each term of the second sum is in $\cI_{S,p}^{\nu_{f,g}^p+1}$, so by definition
$$
\nu_{f,g}^p(h)\geq \text{min}\{\nu_{f,g}^p(\r^1),\dots,\nu_{f,g}^p(\r^k),\nu_{f,g}^p+1\}.
$$
This is true for any $h\in I_{f(p)}$ hence
$$
\text{min}\{\nu_{f,g}^p(h),\text{ for }h\in I_{f(p)}\}\geq \text{min}\{\nu_{f,g}^p(\ell^1),\dots,\nu_{f,g}^p(\ell^k),\nu_{f,g}^p+1\}.
$$
Reminding Definition \ref{def_tang} the assertion follows easily.\par
ii) By Remark \ref{remark1} if $\{\r^1,\dots,\r^k\}$ are generators of $I_{f(p)}$ then the corresponding germs are generators of $I_{f(q)}$ for all $q\in S$ close enough to $p$. Then $ii)$ of Lemma \ref{lemma1} and $i)$ of the current Lemma imply that the set of points of $S$ where $f$ and $g$ are tangential is open and closed at the same time and the assertion follows because $S$ is connected.
\end{pf}
By Lemma \ref{lemma_tang} we can say that {\it $(f,g)$ is tangential along $S$} if it is tangential at some point $p\in S$.\par
Now let assume in the following $S$ non-singular for simplicity and let $p$ be any point of $S$, $(z^1,\dots,z^n)$ any local coordinates adapted to $(g,S)$ at $p$ and $(w^1,\dots,w^n)$ the corresponding special coordinates at $f(p)$. The ideal $I_{f(p)}$ is clearly generated by the germ of $w^1$ at $f(p)$, hence by Lemma \ref{lemma_tang} it easily follows that:
\begin{align*}
\text{$(f,g)$ tangential along $S$ }\Longleftrightarrow &\phantom{a} h^1\in\cI_{S,p},\text{ for every local coordinates}\\
 & \text{adapted to $(g,S)$  at $p$}, \text{ for every }p\in S \\
\Longleftrightarrow & \phantom{a} \sD_{f,g} \text{ is in fact ``tangential to $S$'' i.e. } \\
 & \text{ a section of the bundle of } (N_S^{\ot\nu})^*\ot TS,
\end{align*}
where the $h^1$ are the ones in (\ref{germs_h}). Therefore when $(f,g)$ is tangential along $S$ the canonical section $\sD_{f,g}$ is in fact a $1$-dimensional holomorphic foliation on $S$,
$$
\sD_{f,g}:N_S^{\ot\nu}\longrightarrow TS,
$$ 
and $N_S^{\ot\nu}$ is the tangent bundle of the foliation.  Observe that the possibly singular set of  $\sD_{f,g}$ is clearly $\text{Sing}(f,g)$. We can also think the foliation to be the distribution
$$
\X_{f,g}=\sD_{f,g}\left(N_S^{\ot\nu}\right)\subset TS,
$$
and if $S^0=S-\text{Sing}(f,g)$ then $\X_{f,g}|_{S^0}$ is a line sub-bundle of $TS^0$. If $(z^1,\dots,z^n)$ are any local coordinates adapted to $(g,S)$ defined on an open subset $U\subset M$ such that $U\cap S\neq\emptyset$ then 
\begin{equation}\label{locgen_tang}
X_{f,g}=\sD_{f,g}\left((\partial z^1)^{\nu}\right)=\sum_{j=2}^n \left.h^j\right|_S\left.\frac{\partial}{\partial z^j}\right|_S
\end{equation}
is a local generator of $\sD_{f,g}$ (eq. $\X_{f,g}$) on $U\cap S$, where the $h^j$ are the ones appearing in (\ref{germs_h}). We call it a {\it canonical local generator} of $\sD_{f,g}$ (eq. $\X_{f,g}$).
\begin{remark}
If we consider the couple $(f,\text{Id}_M)$ then clearly it is tangential along $S$ if and only if $f$ is tangential along $S$ according to \cite[Def.1.2.]{ABT2004}. Moreover our foliation $\X_{f,\text{Id}_M}$ is exactly the foliation $\X_f$ of \cite[Def.3.2.]{ABT2004} (when $f$ is tangential along $S$), which they call `canonical distribution'.
\end{remark}
\begin{remark}\label{weaker_hp}
When $(f,g)$ is tangential along $S$ we would be able to define a $1$-dimensional holomorphic foliation on $S$ even weakening the hypothesis on $g$. If we only assume that $f|_S=g|_S:S\to M$ is a local holomorphic embedding, that is
$$
\text{d}f|_{T_pS}=\text{d}g|_{T_pS}:T_pS\longrightarrow T_{f(p)}M
$$
is injective for every $p\in S$, we are able to define a sort of ``pre-canonical section''
$$
\sD_{f,g}^{pre}:N_S^{\ot\nu}\longrightarrow\left.\left(f^*TM\right)\right|_S.
$$
By the hypothesis we have the injective morphism 
$$
D_S=\text{d}f|_{TS}=\text{d}g|_{TS}:TS\longrightarrow\left.\left(f^*TM\right)\right|_S
$$
and one can prove that $(f,g)$ is tangential along $S$ (Remark \ref{remark1} can be suitably adapted and then Lemma \ref{lemma_tang} is still true) if and only if $\sD_{f,g}^{pre}(N_S^{\ot\nu})\subset D_S(TS)$. It follows that in this case we can define the $1$-dimensional holomorphic foliation on $S$
$$
\sD_{f,g}=D_S^{-1}\circ\sD_{f,g}^{pre}:N_S^{\ot\nu}\longrightarrow TS.
$$
On the contrary, if $(f,g)$ is not tangential along $S$ we can not even define a morphism with image into $TM|_S$ (that is the canonical section) and this is one of the reasons we assume the stronger hypothesis on $g$ (the other is that anyway with the weaker hypothesis one can not extend locally the foliation $\sD_{f,g}$ arising in the tangential case - see Section \ref{loc_ext} to understand what we mean).
\end{remark}
If $(f,g)$ is not tangential along $S$ but $S$ {\it sits into $M$ in a particularly nice way} we still have $1$-dimensional holomorphic foliations on $S$. For details about this property of $S$ see \cite[Sec.2]{ABT2004} (or for a deeper treatment \cite[Sec.2 and 3]{ABT2008}), here we just recall some facts and remarks which we are going to use in this paper.\par
\begin{defin}
Let $M$ be a $n$-dimensional complex manifold and $S\subset M$ a complex sub-manifold of any codimension $k$ ($0<k<n$). We say that {\it $S$ splits into $M$} if the short exact sequence
$$
0\longrightarrow TS\stackrel{\imath}{\longrightarrow} \left.TM\right|_S \stackrel{\pi}{\longrightarrow} N_S \longrightarrow 0
$$
splits (holomorphically), that is if there exists a morphism $\s:N_S \to TM|_S$ of holomorphic vector bundles such that $\pi\circ\s=\text{Id}_{N_S}$.\par
Let $\gU$ be an atlas of $M$ adapted to $S$. We call it a {\it splitting atlas adapted to $S$} if  
\begin{equation}\label{split_prop}
\frac{\partial \hat{z}^i}{\partial z^j}\in\cI_S(U\cap\hat{U}),\qquad\qquad i=k+1,\dots,n,\phantom{a}j=1,\dots,k,
\end{equation}
for any two charts $(U,z)$ and $(\hat{U},\hat{z})$ in $\gU$ such that $U\cap\hat{U}\neq\emptyset$.
\end{defin}
Then
\begin{proposition}
Let $M$ be a $n$-dimensional complex manifold and $S\subset M$ a complex sub-manifold of any codimension $k$ ($0<k<n$). The following statements are equivalent:
\begin{itemize}
\item[(i)] $S$ splits into $M$
\item[(ii)] there exists a morphism $\t:TM|_S\to TS$ of holomorphic vector bundles such that $\t\circ\imath=\text{Id}_{TS}$
\item[(iii)] there exists a splitting atlas adapted to $S$
\end{itemize}
\end{proposition}
\begin{pf}
(i)$\Longleftrightarrow$ (iii): see \cite[Prop.2.1.]{ABT2004} or even \cite[Prop.2.15.]{ABT2008}.\\
(i)$\Longleftrightarrow$ (ii): suppose there exists $\s$, then $\pi\circ (\text{Id}_{TM|_S}-\s\circ\pi)=0$ by the property of $\s$ and this implies that $\text{im}(\text{Id}_{TM|_S}-\s\circ\pi)\subset\text{ker}(\pi)=\text{im}(\imath)$. Hence $\t=\imath^{-1}\circ (\text{Id}_{TM|_S}-\s\circ\pi)$ does the work.\par
Conversely, suppose there exists $\t$. By its property it follows that $\text{ker}(\t)\cap\text{im}(\imath)=\text{ker}(\t)\cap\text{ker}(\pi)=\{\text{zero section}\}$ hence we can invert $\pi|_{\text{ker}(\t)}$ and $\s=\pi|_{\text{ker}(\t)}^{-1}$.
\end{pf}
Observe that by definitions $\t\circ\s=\imath^{-1}\circ (\s-\s)=0$, hence $\text{im}(\s)\subset\text{ker}(\t)$. Moreover $\text{ker}(\t)=\text{ker}(\text{Id}_{TM|_S}-\s\circ\pi)$ then $\text{im}(\s)=\text{ker}(\t)$.\par
\begin{example}
Let $S\subset M$ be a sub-manifold and $U\subset M$ an open neighborhood of $S$ in $M$. If there exists a holomorphic retraction $\r:U\to S$ then clearly $S$ splits into $M$.
\end{example}
\begin{example}\label{example}
A rank $k$ holomorphic vector bundle $\pi:M\to S$ is a holomorphic retraction of $M$ on $S$ (just identify $S$ with the image of the zero section of the bundle). Thus the base of a holomorphic vector bundle can be seen as a splitting sub-manifold of the total space of the bundle. The local charts on $M$ induced by a trivialization of the bundle clearly give a splitting atlas adapted to $S$.
\end{example}
Let $S\subset M$ be a sub-manifold of any codimension $k$ ($0<k<n$) and $\s:N_S\to TM|_S$ any {\it splitting morphism} of $S$ into $M$. If $\gU$ is any splitting atlas adapted to $S$ and $(U,z)\in\gU$ any local chart then $\{\partial z^1=\pi(\frac{\partial}{\partial z^1}|_S),\dots,\partial z^k=\pi(\frac{\partial}{\partial z^k}|_S)\}$ is a local frame for $N_S$ but in general $\s(\partial z^j)\neq \frac{\partial}{\partial z^j}|_S$ for $j=1,\dots,k$ (we only know that their differences are in $\cT_S$). This observation leads to the following definition.
\begin{defin}\label{s-split-atlas}
Let $M$ be a $n$-dimensional complex manifold and $S\subset M$ a complex sub-manifold of codimension $k$ ($0<k<n$) which splits into $M$. Let $\s$ be a splitting morphism. An atlas $\gU$ of $M$ is said to be a {\it $\s$-splitting atlas adapted to $S$} if it is a splitting atlas adapted to $S$ such that $\s(\partial z^j)= \frac{\partial}{\partial z^j}|_S$ for $j=1,\dots,k$.\par
If moreover there is a local biholomorphism $\f:W\to M$ defined on an open neighborhood $W\subset M$ of $S$ and every $(U,z)\in\gU$ is such that $\left. \f\right|_U$ is a biholomorphism onto its image, then we call $\gU$ a {\it $\s$-splitting atlas adapted to $(\f,S)$}.
\end{defin}
One can easily show that given $S\subset M$ and a splitting morphism $\s:N_S\to TM|_S$ then a {\it $\s$-splitting atlas adapted to $S$} always exists (hence if there is also a  $\f:W\to M$ as in Definition \ref{s-split-atlas} a $\s$-splitting atlas adapted to $(\f,S)$ always exists).\par
Now let $S\subset M$ be of codimension $1$ again and let $(f,g)\in\text{End}_S^2(M)$ be a couple whit $g$ a local biholomorphism around $S$, whose order of coincidence is $\nu=\nu_{f,g}$. If $S$ splits into $M$ then for any splitting morphism $\s$ we can define the $1$-dimensional holomorphic foliation on $S$
$$ 
\sD_{f,g}^{\s}=\t^{\s}\circ\sD_{f,g}:N_S^{\ot\nu}\longrightarrow TS,
$$ 
where $\t^{\s}=\imath^{-1}\circ (\text{Id}_{TM|_S}-\s\circ\pi)$. The tangent bundle of this foliation is $N_S^{\ot\nu}$ and its (possibly) singular set $\text{Sing}(\sD_{f,g}^{\s})$ may be larger than $\text{Sing}(f,g)$. As before, we can also think the foliation to be the distribution
$$
\X_{f,g}^{\s}=\sD_{f,g}^{\s}\left(N_S^{\ot\nu}\right)\subset TS,
$$
and if $S^0=S-\text{Sing}(\sD_{f,g}^{\s})$ then $\X_{f,g}^{\s}|_{S^0}$ is a line sub-bundle of $TS^0$. If $\gU$ is a $\s$-splitting atlas adapted to $(g,S)$ and $(U,z)\in\gU$ such that $U\cap S\neq\emptyset$
$$
\t^{\s}\left(\left.\frac{\partial}{\partial z^1}\right|_S\right)=\t^{\s}\circ\s\left(\partial z^1\right)=0,
$$
then
\begin{equation}\label{locgen_tang_split}
X_{f,g}^{\s}=\sD_{f,g}^{\s}\left((\partial z^1)^{\nu}\right)=\sum_{j=2}^n \left.h^j\right|_S\left.\frac{\partial}{\partial z^j}\right|_S
\end{equation}
is a local generator of $\sD_{f,g}^{\s}$ (eq. $\X_{f,g}^{\s}$) on $U\cap S$, where the $h^j$ are the ones appearing in (\ref{germs_h}). We call it a {\it canonical local generator} of $\sD_{f,g}^{\s}$ (eq. $\X_{f,g}^{\s}$).\par
When $\nu_{f,g}=1$ we can define even other $1$-dimensional holomorphic foliations on $S$. By the hypothesis on $(f,g)$ we have the morphisms of holomorphic vector bundles on $S$ 
$$
\text{d}f|_S,\text{ d}g|_S:\left. TM\right|_S\longrightarrow \left. \left(f^*TM\right)\right|_S
$$
which are different but coincide on $TS\subset TM|_S$. Since $\text{d}g|_S$ is invertible we can compose $\text{d}g|_S^{-1}\circ\text{d}f|_S$ and it induces the morphism
\begin{align}\label{dfg}
d_{f,g}:N_S &  \longrightarrow N_S \\
[v] & \longrightarrow [(\text{d}g_p^{-1}\circ\text{d}f_p)(v)],\qquad v\in T_pM,p\in S. \notag
\end{align}
As a consequence we can define the $1$-dimensional holomorphic foliations 
$$
\sD_{f,g}^{\s,1}=\sD_{f,g}^{\s}\circ d_{f,g}:N_S\longrightarrow TS
$$
on $S$, where $\s$ is any splitting morphism of $S$ in $M$ as before.
\begin{remark}\label{rmk_A}
One might consider $d_{f,g}^{\ot\nu}:N_S^{\ot\nu}\to N_S^{\ot\nu}$ for any $\nu\geq 1$ and define the foliations
$$
\sD_{f,g}^{\s,\nu}=\sD_{f,g}^{\s}\circ d_{f,g}^{\ot\nu}:N_S^{\ot\nu}\longrightarrow TS
$$ 
but since $d_{f,g}=\text{Id}_{N_S}$ if and only if $(f,g)$ tangential along $S$ or $\nu>1$ then it would be $\sD_{f,g}^{\s,\nu}=\sD_{f,g}^{\s}$ for any $\nu>1$.
\end{remark}
The tangent bundle of this foliation is $N_S$ and its (possibly) singular set $\text{Sing}(\sD_{f,g}^{\s,1})$ may be larger than $\text{Sing}(\sD_{f,g}^{\s})$. As in the other cases we can also think the foliation to be the distribution
$$
\X_{f,g}^{\s,1}=\sD_{f,g}^{\s,1}\left(N_S\right)\subset TS
$$
and if $S^0=S-\text{Sing}(\sD_{f,g}^{\s,1})$ then $\X_{f,g}^{\s,1}|_{S^0}$ is a line sub-bundle of $TS^0$. If $\gU$ is a $\s$-splitting atlas adapted to $(g,S)$ and $(U,z)\in\gU$ such that $U\cap S\neq\emptyset$ then one can easily check that
$$
d_{f,g}\left(\partial z^1\right)=\left.\frac{\partial f^1}{\partial z^1}\right|_S\partial z^1=\left(1+h^1|_S\right)\partial z^1
$$
where $f^1=w^1\circ f$, with $(w^1,\dots,w^n)$ special coordinates induced by the coordinates $z$ and $h^1$ is the one appearing in (\ref{germs_h}). Hence
\begin{equation}\label{locgen_tang_split_nu1} 
X_{f,g}^{\s,1}=\sD_{f,g}^{\s}\left(\left(1+h^1|_S\right)\partial z^1\right)=\left(1+h^1|_S\right)\sum_{j=2}^n \left.h^j\right|_S\left.\frac{\partial}{\partial z^j}\right|_S
\end{equation}
is a local generator of $\sD_{f,g}^{\s,1}$ (eq. $\X_{f,g}^{\s,1}$) on $U\cap S$, where the $h^j$ are the ones appearing in (\ref{germs_h}). We call it a {\it canonical local generator} of $\sD_{f,g}^{\s,1}$ (eq. $\X_{f,g}^{\s,1}$).\par
\begin{remark}
Suppose that $S$ splits into $M$ and $(f,g)$ is tangential along $S$ at the same time and let $\s$ be any splitting morphism. Since $\t^{\s}\circ\imath=\text{Id}_{TS}$ it follows that $\sD_{f,g}^{\s}=\sD_{f,g}$ for any $\s$ and $\nu\geq 1$. Moreover, if $\nu=1$ we have also $\sD_{f,g}^{\s,1}=\sD_{f,g}$ by Remark \ref{rmk_A}.
\end{remark}
We conclude by coming back to Example \ref{example} for a while, in order to introduce a definition which we will use from Section \ref{loc_ext} on.\par
If $\pi:M\to S$ is a rank $k$ holomorphic vector bundle then $S$ is something more of a splitting sub-manifold of $M$. In fact, if $(U,z)$ and $(\hat{U},\hat{z})$ are two local charts on $M$ induced by a trivialization of the bundle, where $(z^{k+1},\dots,z^n)$ and $(\hat{z}^{k+1},\dots,\hat{z}^n)$ are local coordinates respectively on $\pi(U)\subset S$ and $\pi(\hat{U})\subset S$, then 
$$
\frac{\partial \hat{z}^i}{\partial z^j}\equiv 0 \qquad\quad\text{ and }\qquad\quad \frac{\partial^2 \hat{z}^r}{\partial z^s\partial z^t}\equiv 0
$$
for $i=k+1,\dots,n$, $j=1,\dots,k$ and $r,s,t=1,\dots,k$. This observation leads to the following definitions.
\begin{defin}
Let $M$ be a $n$-dimensional complex manifold and $S\subset M$ a complex sub-manifold of codimension $k$ ($0<k<n$). We say that {\it $S$ is comfortably embedded into $M$} if there exists a splitting atlas $\gU$ adapted to $S$ (hence $S$ splits into $M$) such that 
\begin{equation}\label{comf_prop}
\frac{\partial^2 \hat{z}^r}{\partial z^s\partial z^t}\in\cI_S(U\cap\hat{U}),\qquad\qquad r,s,t=1,\dots,k,
\end{equation}
for any two charts $(U,z)$ and $(\hat{U},\hat{z})$ in $\gU$ such that $U\cap\hat{U}\neq\emptyset$.\par
Such an atlas is said to be a {\it comfortably atlas adapted to $S$}. If moreover $\s$ is a splitting morphism of $S$ in $M$ and $\gU$ is also a $\s$-splitting atlas adapted to $S$ then it is called a {\it $\s$-comfortably atlas adapted to $S$}. Finally, if there is a local biholomorphism $\f:W\to M$ defined on an open neighborhood $W\subset M$ of $S$ and every $(U,z)\in\gU$ is such that $\left. \f\right|_U$ is a biholomorphism onto its image, then we call $\gU$ a {\it $\s$-comfortably atlas adapted to $(\f,S)$}.
\end{defin}
Roughly speaking, a comfortably embedded sub-manifold is a sort of first-order approximation of the zero section of a holomorphic vector bundle.\par
\begin{example}
Let $M$ be a $n$-dimensional complex manifold and $p\in M$ some point. Let $\pi:\tilde{M}\to M$ denote the blow-up of $M$ at $p$ and let $S=\pi^{-1}(p)\cong\bP^{n-1}$ be the {\it exceptional divisor}, which is a non-singular compact connected hypersurface in $\tilde{M}$. Then it is easy to check that $S$ is comfortably embedded into $\tilde{M}$.
\end{example}
\begin{remark}
In both Section \ref{ord_can_sec} and \ref{foliations} we have assumed $S\subset M$ to be a hypersurface, except for some general definitions. Anyway all definitions, lemmas, propositions and computations stated and done up to now may be easily adapted
for $S$ of any codimension $k$ ($0<k<n$), likewise in the first three sections of \cite{ABT2004}. Conversely, we need $S$ to have codimension $1$ for what follows.
\end{remark}
\bigskip
\section{A Baum-Bott-type index theorem}\label{BaumBott}\setcounter{equation}{0}
Let $S$ be a $m$-dimensional complex manifold and suppose to have a $1$-dimensional holomorphic (possibly singular) foliation 
$$
\sF:F\to TS
$$ 
on it. Set $S^0=S-\text{Sing}(\sF)$ and $F^0=F|_{S^0}$, which we can identify with the line sub-bundle $\sF(F^0)\subset TS^0$.
\begin{defin}\label{normal_fol}
The quotient $N_{\sF}=TS^0/F^0$ is called {\it normal bundle of $\sF$}, while the virtual bundle (in the sense of $K$-theory) $TS-F$ is the {\it virtual normal bundle of $\sF$}, since $(TS-F)|_{S^0}=N_{\sF}$ in the $K$-group of $S^0$.\par
If we express the foliation through the injective morphism $\sF:\cF\to\cT_S$ then the coherent $\cO_S$-module  $\cN_{\sF}=\cT_S/\cF$ is called {\it normal sheaf of $\sF$}. Clearly $\cF^0=\cF|_{S^0}$ is the sheaf of germs of holomorphic sections of $F^0$ and $\cN_{\sF}^0=\cN_{\sF}|_{S^0}$ is the sheaf of germs of holomorphic sections of $N_{\sF}$.
\end{defin}
Let $\pi:TS^0\to N_{\sF}$ be the obvious projection. There is the natural partial holomorphic connection (in the sense of Bott \cite{bo1967}) on $N_{\sF}$ along $F^0$ 
\begin{align}\label{BB_conn}
\cN_{\sF}^0 & \stackrel{\d^{bb}}{\longrightarrow} \left(\cF^0\right)^*\ot\cN_{\sF}^0 \notag \\
w & \longrightarrow \d^{bb}(w)\text{ s.t. }\d^{bb}(w)(v)=\pi\left([v,\tilde{w}]\right), 
\end{align}
for any $w\in\cN_{\sF}^0$ and $v\in\cF^0$, where $\tilde{w}\in\cT_{S^0}$ is any vector field such that $\pi(\tilde{w})=w$. Observe that $\d^{bb}(w)(v)$ is holomorphic whenever $w$ and $v$ are, and (\ref{BB_conn}) is independent by the choice of $\tilde{w}$.\par
This partial connection makes $N_{\sF}$ a {\it $F^0$-bundle} (using the terminology of \cite{bo1967}, see also \cite[Sec.II.9.]{Su1998}) and as a consequence of the `Bott vanishing theorem' (see \cite[Th.II.9.11.]{Su1998} for a complete version or \cite[Th.6.2.3.]{BSS2009} for a simplified one) and of the the exact sequence
$$
0\longrightarrow F^0\stackrel{\sF|_{S^0}}{\longrightarrow} TS^0 \stackrel{\pi}{\longrightarrow} N_{\sF}\longrightarrow 0,
$$
one can localize at $\text{Sing}(\sF)$ suitable characteristic classes of the virtual normal bundle $TS-F$. In particular if $S$ is {\it compact} one gets theorem \cite[Th.1.]{BB1972} (see also \cite[Th.III.7.6.]{Su1998} or \cite[Th.6.2.5.]{BSS2009}). Then by it and by Section \ref{foliations} we gain the following index theorem.
\begin{theorem}[Baum-Bott-type index theorem]\label{BB_theorem}
Let $M$ be a $n$-dimensional complex manifold and $S\subset M$ a non-singular compact connected complex hypersurface in $M$. Let $(f,g)\in\text{End}_S^2(M)$ be a couple where $g$ is a local biholomorphism around $S$ and set $\nu=\nu_{f,g}$. Assume that
\begin{itemize}
\item[(i)] $(f,g)$ is tangential along $S$ 
\end{itemize}
or that
\begin{itemize}
\item[(ii)] $S$ splits into $M$.
\end{itemize}
In case $(i)$ set $\sD=\sD_{f,g}$ while in case $(ii)$ set $\sD=\sD_{f,g}^{\s}$ (where $\s$ is some splitting morphism) and suppose moreover that $\sD\neq 0$ and $\sD^{\s,1}_{f,g}\neq 0$. Let $\text{Sing}(\sD)=\sqcup_{\l}\S_{\l}$ and $\text{Sing}(\sD^{\s,1}_{f,g})=\sqcup_{\mu}\S_{\mu}^1$  be the decompositions in connected components of the singular sets of the foliations. \par
Then for any symmetric homogeneous polynomial $\f\in\bC[z_1,\dots,z_{n-2}]$ of degree $n-1$ there exist complex numbers $\text{Res}_{\f}(\sD;TS-N_S^{\ot\nu};\S_{\l})$ and $\text{Res}_{\f}(\sD^{\s,1}_{f,g};TS-N_S;\S_{\mu}^1)$ such that
$$
\sum_{\l}\text{Res}_{\f}\left(\sD;TS-N_S^{\ot\nu};\S_{\l}\right)=\int_S\f\left(TS-N_S^{\ot\nu}\right)
$$
and, only when $\nu=1$,
$$
\sum_{\mu}\text{Res}_{\f}\left(\sD^{\s,1}_{f,g};TS-N_S;\S_{\mu}^1\right)=\int_S\f\left(TS-N_S\right).
$$
\end{theorem}   
We call Theorem \ref{BB_theorem} a {\it Baum-Bott-type index theorem} because Baum and Bott have introduced this kind of residues and the partial connection (\ref{BB_conn}) (see \cite{BB1972} and \cite{BB1970}). Moreover, as we said, the theorem follows by the Baum-Bott index theorem \cite[Th.1.]{BB1972}.
\begin{remark}
To be precise in case $(ii)$ one just needs that $S-\text{Sing}(f,g)$ splits into $M$ since $\text{Sing}(f,g)\subset\text{Sing}(\sD_{f,g}^{\s})$.
\end{remark}
\begin{remark}
We can think the foliations of Theorem \ref{BB_theorem} in terms of morphisms of $\cO_S$-modules, that is as
$$
\sD:\cN_S^{\ot\nu}\to\cT_S\qquad\text{ and }\qquad \sD^{\s,1}_{f,g}:\cN_S\to\cT_S,
$$
which are injective morphisms. Recall that the coherent $\cO_S$-modules $\cN_{\sD}=\cT_S/\sD(\cN_S^{\ot\nu})$ and $\cN_{\sD}^1=\cT_S/\sD^{\s,1}_{f,g}(\cN_S)$ are the normal sheaves of $\sD$ and $\sD^{\s,1}_{f,g}$ and their restrictions to $S^0$ are the sheaves of germs of holomorphic sections respectively of $N_{\sD}$ and $N_{\sD^{\s,1}_{f,g}}$. Then we have 
$$\f\left(TS-N_S^{\ot\nu}\right)=\f\left(\cN_{\sD}\right) \qquad\text{ and }\qquad \f\left(TS-N_S\right)=\f\left(\cN_{\sD}^1\right)
$$
for any symmetric homogeneous polynomial $\f$ (see \cite[Ch.VI]{Su1998}).
\end{remark}
\begin{remark}
If we consider the couple $(f,\text{Id}_M)$ then Theorem \ref{BB_theorem} turns out to be \cite[Th.6.4.]{ABT2004}. Observe that Abate-Bracci-Tovena assume $S-\text{Sing}(f,\text{Id}_M)$ to be comfortably embedded into $M$ and not only splitting, however the splitting property is sufficient. Moreover they do not consider the foliation $\sD_{f,\text{Id}}^{\s}$ when $\nu=1$ but only $\sD_{f,\text{Id}}^{\s,1}$ (which in their notation are respectively $H_{\s,f}$ and $H_{\s,f}^1$).
\end{remark}
\begin{remark}
By Remark \ref{weaker_hp} when $(f,g)$ is tangential along $S$ we may assume only that $f|_S=g|_S:S\to M$ is a local holomorphic embedding (instead of the stronger ``$g$ is a local biholomorphism around $S$'') and get anyway an index theorem likewise Theorem \ref{BB_theorem} yet.
\end{remark}
We conclude this section by deriving explicit formulas for the computation of the residues in Theorem \ref{BB_theorem} at isolated singular points. For this purpose we briefly recall how these residues are defined in Lehmann-Suwa theory. We do it considering only foliations $\sD$ for simplicity, anyway for $\sD^{\s,1}_{f,g}$ it is the same. As a reference consult \cite[Sec.III.7.]{Su1998}.\par
Let $\f\in\bC[z_1,\dots,z_{n-2}]$ be any symmetric homogeneous polynomial of degree $n-1$ and $\d_{\sD}^{bb}$ be the partial holomorphic connection on $N_{\sD}$ along $N_{S^0}^{\ot\nu}$ defined as in (\ref{BB_conn}). Let $\n_{\sD}^{bb}$ be any $(1,0)$-type extension of it, which always exists (see for example \cite[Lem.2.5.]{BB1972}). Then put connections $\n_1^0$ and $\n_2^0$ respectively on $N_{S^0}^{\ot\nu}$ and $TS^0$ such that the triple $(\n_1^0,\n_2^0,\n_{\sD}^{bb})$ is compatible (see \cite[p.72]{Su1998} for the meaning) with the short exact sequence
\begin{equation}\label{BB_seq}
0\longrightarrow N_{S^0}^{\ot\nu}\stackrel{\sD|_{S^0}}{\longrightarrow}TS^0\stackrel{\pi}{\longrightarrow} N_{\sD}\longrightarrow 0
\end{equation}
and set $\n_{\sD}^{\bullet}=(\n_1^0,\n_2^0)$. One makes this choices since by the compatibility of the triple and the Bott vanishing theorem
$$
\f\left(\n_{\sD}^{\bullet}\right)=\f\left(\n_{\sD}^{bb}\right)=0,
$$
where $\f(\n_{\sD}^{\bullet})$ is a $2(n-1)$ form on $S^0$ defined as at \cite[pp.71-72]{Su1998}. Let now $\S$ be a connected component of $\text{Sing}(\sD)$ and $V\subset S$ an open  sub-set such that $V\cap \text{Sing}(\sD)=\S$. Choose any connections $\n_1^V$ and $\n_2^V$ on $N_S^{\ot\nu}|_V$ and $TS|_V$ and set $\n_V^{\bullet}=(\n_1^V,\n_2^V)$. Finally, let $R\subset V$ be a compact real sub-manifold of dimension $2(n-1)$ oriented as $S$ and such that $\S\subset\text{int}(R)$. Consider on the boundary $\partial R$ the orientation induced by $R$. Then by definition the residue is
$$
\text{Res}_{\f}\left(\sD;TS-N_S^{\ot\nu};\S\right)=\int_R\f\left(\n_V^{\bullet}\right)-\int_{\partial R}\f\left(\n_{\sD}^{\bullet},\n_V^{\bullet}\right),
$$
where the Bott difference form $\f(\n_{\sD}^{\bullet},\n_V^{\bullet})$ is a $2(n-1)-1$ form on $V-\S$ defined as at \cite[pp.71-72]{Su1998}. One can show that this formula does not depend on the choice of the various connections or of the sub-manifold $R$.\par
If $\S=\{p\}$ is an isolated point we can assume $V$ to be such that $N_S^{\ot\nu}|_V$ and $TS|_V$ are trivial. Hence if we take $\n_1^V$ and $\n_2^V$ trivial respect to some local frames $\f(\n_V^{\bullet})= 0$ and the residue becomes
$$
\text{Res}_{\f}\left(\sD;TS-N_S^{\ot\nu};p\right)=-\int_{\partial R}\f\left(\n_{\sD}^{\bullet},\n_V^{\bullet}\right).
$$
Observe now that $N_{S^0}^{\ot\nu}$ and $TS^0$ are not in general $N_{S^0}^{\ot\nu}$-bundles but if there is a local generator $X$ of the foliation $\sD$ on $V$ then 
$$
\left.N_{S^0}^{\ot\nu}\right|_{V-\{p\}}\qquad\text{ and }\qquad \left.TS^0\right|_{V-\{p\}}
$$
are canonically $N_{S^0}^{\ot\nu}|_{V-\{p\}}$-bundles thanks to the natural `holomorphic action'  of $X$ on them by Lie bracket $[X,\cdot]$ (again we are using the terminology of \cite[Sec.II.9.]{Su1998}). This action induces partial holomorphic connections on $N_{S^0}^{\ot\nu}|_{V-\{p\}}$ and $TS^0|_{V-\{p\}}$ along $N_{S^0}^{\ot\nu}|_{V-\{p\}}$ which are, together with $\d_{\sD}^{bb}$, compatible with (\ref{BB_seq}) restricted to $V-\{p\}$. Therefore we can assume $\n_1^0$ and $\n_2^0$ to be `$X$-connections' defined on $V-\{p\}$ (and again such that the triple $(\n_1^0,\n_2^0,\n_{\sD}^{bb})$ is compatible with (\ref{BB_seq})). With this choice one can show that
$$
\text{Res}_{\f}\left(\sD;TS-N_S^{\ot\nu};p\right)=-\int_{\partial R}\f\left(\n_2^0,\n_2^V\right)
$$
and using the same arguments for $\sD^{\s,1}_{f,g}$ also that
$$
\text{Res}_{\f}\left(\sD^{\s,1}_{f,g};TS-N_S;p\right)=-\int_{\partial R}\f\left(\n_2^0,\n_2^V\right),
$$
where $\n_2^0$ and $\n_2^V$ are as above.\par
Now one can work as in the proof of \cite[Th.III.5.5.]{Su1998} (see also \cite[Sec.5.]{le1991}) and obtain a similar formula. In particular, let $(U,z)$ be a local chart of $M$ at $p$ belonging to an atlas $\gU$ adapted to $(g,S)$ (or a $\s$-splitting atlas adapted to $(g,S)$, if necessary), set $V=U\cap S$ and let $X$ be the canonical local generators (\ref{locgen_tang}), (\ref{locgen_tang_split}) or (\ref{locgen_tang_split_nu1}) (depending on the case). Moreover let $\n_2^V$ be trivial respect to the local frame $\{\frac{\partial}{\partial z^2}|_V,\dots,\frac{\partial}{\partial z^n}|_V\}$ of $TS$ and let the $h^j\in\cO_{M,p}$ be the ones appearing in (\ref{germs_h}). Then when $(f,g)$ is tangential along $S$ and $\sD=\sD_{f,g}$  or when $S$ splits into $M$ and $\sD=\sD_{f,g}^{\s}$ we get the formula
\begin{equation}\label{BB_formula_1}
\text{Res}_{\f}\left(\sD;TS-N_S^{\ot\nu};p\right)=\int_{\G}\frac{\f\left(-H\right)}{h^2\cdots h^n}\text{ d}z^2\wedge\cdots\wedge\text{d}z^n,
\end{equation}
where 
$$
H=\left(\left. \frac{\partial h^j}{\partial z^k}\right|_{U\cap S}\right)_{j,k=2,\dots, n}
$$
and $\G$ is the $(n-1)$ cycle
$$
\G=\left\{q\in U\cap S\text{ s.t. }\left|h^2(q)\right|=\dots=\left|h^n(q)\right|=\e\right\},
$$
for $\e>0$ small enough, oriented so that $\text{d}\q^2\wedge\cdots\wedge\text{d}\q^n>0$ where $\q^j=\text{arg}(h^j)$. Similarly, when $S$ splits into $M$ and $\nu=1$ we have
\begin{equation}\label{BB_formula_2}
\text{Res}_{\f}\left(\sD^{\s,1}_{f,g};TS-N_S;p\right)=\int_{\G'}\frac{\f\left(-H'\right)}{\left(1+h^1\right)^{n-1}h^2\cdots h^n}\text{ d}z^2\wedge\cdots\wedge\text{d}z^n,
\end{equation}
where 
$$
H'=\left(\left. \frac{\partial \left(1+h^1\right)h^j}{\partial z^k}\right|_{U\cap S}\right)_{j,k=2,\dots, n}
$$
and
$$
\G'=\left\{q\in U\cap S\text{ s.t. }\left|(1+h^1(q))h^2(q)\right|=\dots=\left|(1+h^1(q))h^n(q)\right|=\e\right\},
$$
for $\e>0$ small enough, oriented so that $\text{d}\q^2\wedge\cdots\wedge\text{d}\q^n>0$ where $\q^j=\text{arg}(1+h^1)(h^j)$.
\begin{remark}
Implicitly we have proved (and used) that
$$
\text{Res}_{\f}\left(\sD;TS-N_S^{\ot\nu};p\right)=\text{Res}_{\f}\left(X;TS;p\right) 
$$
and
$$ 
\text{Res}_{\f}\left(\sD^{\s,1}_{f,g};TS-N_S;p\right)=\text{Res}_{\f}\left(X;TS;p\right),
$$
where $X$ is a local generator of the foliations $\sD$ or $\sD^{\s,1}_{f,g}$ and the residues on the right are the ones associated to the natural `holomorphic action' of $X$ on $TS^0|_{V-\{p\}}$ by Lie bracket (see \cite[Sec.III.5.]{Su1998} for definition). Consequently we could have used directly the formula of \cite[Th.III.5.5.]{Su1998}. See also \cite[Rmk.III.7.7.(1)]{Su1998} or \cite[Rmk.6.2.2(1)]{BSS2009}.
\end{remark}
\bigskip
\section{Local extensions of the foliations}\label{loc_ext}\setcounter{equation}{0}
Let $M$ be a $n$-dimensional complex manifold and $S\subset M$ an analytic sub-variety of pure dimension $m$ with regular part $S'=S-\text{Sing}(S)$. Suppose to have a $1$-dimensional holomorphic (possibly singular) foliation on $S'$
$$
\sF:F\to TS'
$$ 
leaving an extension to an open neighborhood $W\subset M$ of $S$, that is there exists a $1$-dimensional holomorphic foliation
$$
\tilde{\sF}:\tilde{F}\to \left.TM\right|_W
$$
on $W$ such that $\tilde{\sF}|_{S'}\equiv\sF$ (hence {\it $\tilde{\sF}$ leaves $S$ invariant}). Let denote $S^0=S'-\text{Sing}(\sF)$ and $F^0=F|_{S^0}$ (which again is identified with $\sF(F^0)\subset TS^0$).
\begin{defin}\label{normal_fol_M}
The quotient $N_{\sF}^M=\left.TM\right|_{S^0}/F^0$ is called {\it normal bundle of $\sF$ in $M$}. Observe that $N_{\sF}^M=N_{\tilde{\sF}}|_{S^0}$ where $N_{\tilde{\sF}}$ is the normal bundle of $\tilde{\sF}$ as in Definition \ref{normal_fol}, defined on $W-\text{Sing}(\tilde{\sF})$.
\end{defin}
In this situation there are natural partial holomorphic connections along $F^0$ on the normal bundle of $S^0$ in $M$, namely $N_{S^0}=\left.TM\right|_{S^0}/TS^0$, and on $N_{\sF}^M$. Indeed, if $\pi:TM|_{S^0}\to N_{S^0}$ is the obvious projection one can define 
\begin{align}\label{CS_conn}
\cN_{S^0} & \stackrel{\d^{cs}}{\longrightarrow} \left(\cF^0\right)^*\ot\cN_{S^0} \notag \\
w & \longrightarrow \d^{cs}(w)\text{ s.t. }\d^{cs}(w)(v)=\pi\left(\left.[\tilde{v},\tilde{w}]\right|_{S^{0}}\right)
\end{align}
for any $w\in\cN_{S^0}$ and $v\in\cF^0$, where $\tilde{w}\in\cT_M|_{S^0}$ and $\tilde{v}\in\tilde{\cF}|_{S^0}$ are any vector fields such that $\pi(\tilde{w}|_{S^{0}})=w$ and $\tilde{v}|_{S^{0}}=v$. Observe that $\d^{cs}(w)(v)$ is holomorphic whenever $w$ and $v$ are and (\ref{CS_conn}) is independent by the choices of $\tilde{w}$ and $\tilde{v}$. Similarly, if $\r:TM|_{S^0}\to N_{\sF}^M$ is the other obvious projection and $\cO(-)$ denotes the sheaves of germs of holomorphic sections of bundles, one can define 
\begin{align}\label{LS_conn}
\cO\left(N_{\sF}^M\right) & \stackrel{\d^{ls}}{\longrightarrow} \left(\cF^0\right)^*\ot\cO\left(N_{\sF}^M\right) \notag \\
w & \longrightarrow \d^{ls}(w)\text{ s.t. }\d^{ls}(w)(v)=\r\left(\left.[\tilde{v},\tilde{w}]\right|_{S^{0}}\right)
\end{align}
for any $w\in\cO(N_{\sF}^M)$ and $v\in\cF^0$, where $\tilde{w}\in\cT_M|_{S^0}$ and $\tilde{v}\in\tilde{\cF}|_{S^0}$ are any sections such that $\r(\tilde{w}|_{S^{0}})=w$ and $\tilde{v}|_{S^{0}}=v$. Again $\d^{ls}(w)(v)$ is holomorphic whenever $w$ and $v$ are and (\ref{LS_conn}) is independent by the choices of $\tilde{w}$ and $\tilde{v}$.\par
These partial connections, possibly with additional hypotheses on $S$, allow one to localize at $\text{Sing}(S)\cup\text{Sing}(\sF)$ suitable characteristic classes of some holomorphic vector bundles defined in a neighborhood of $S$ and then to get index theorems when $S$ is compact. (see \cite[Sec.IV.5 and IV.6.]{Su1998} or \cite[Sec.6.3.2 and 6.3.3]{BSS2009} for details).\par
\begin{remark}
In the setting just described the double inclusion $F^0\subset TS^0 \subset TM|_{S^0}$ induces the ``normal exact sequence''
$$
0\longrightarrow N_{\sF} \longrightarrow N_{\sF}^M \longrightarrow N_{S^0} \longrightarrow 0
$$
and the canonical partial holomorphic connections along $F^0$ defined in (\ref{BB_conn}), (\ref{CS_conn}) and (\ref{LS_conn}) are compatible with it (see \cite[p.72]{Su1998} for the meaning). 
\end{remark}
Let $S$ be a globally irreducible hypersurface and $(f,g)\in\text{End}_S^2(M)$ a couple where $g$ is a local biholomorphism around $S'$. If $(f,g)$ is tangential along $S$ or if $S'$ splits into $M$ then we know by Section \ref{foliations} that there are natural foliations on $S'$, namely $\sD_{f,g}$, $\sD_{f,g}^{\s}$ and $\sD_{f,g}^{\s,1}$. In general they cannot be extended to foliations on a whole neighborhood of $S$ but we can do something weaker. Indeed we can extend them {\it locally} (about points of $S'$) in a quite good way, that is so that we are able to define partial holomorphic connections almost as in (\ref{CS_conn}) and (\ref{LS_conn}). We are going to talk about this sort of ``first order extensions'' of foliations here, while the resulting partial holomorphic connections will be discussed in Section \ref{CamachoSad} and \ref{LehmannSuwa} (see \cite{br2004} to investigate further this topic). \par
\smallskip
For simplicity, in the following let $S$ be a non-singular hypersurface. Assume $(f,g)$ to be tangential along $S$ with order of coincidence $\nu=\nu_{f,g}$ and fix an atlas $\gU$ on $M$ adapted to $(g,S)$. Recall that if $(U,z)\in\gU$ is a local chart such that $U\cap S\neq\emptyset$ then (\ref{locgen_tang}) gives the canonical local generator $X_{f,g}$ of the foliation $\sD_{f,g}$ (eq. $\X_{f,g}$) associated to these coordinates. We now define the local holomorphic vector field
\begin{equation}\label{can_gen_tang_ext}
\cX_{f,g}=\gD_{f,g}\left(\left(\frac{\partial}{\partial z^1}\right)^{\ot\nu}\right)=\sum_{j=1}^n h^j \frac{\partial}{\partial z^j},
\end{equation}
where $\gD_{f,g}$ is defined in (\ref{XXX}), to be the {\it canonical local extension} of $X_{f,g}$. Observe that it generates a $1$-dimensional holomorphic foliation on $U$ leaving $U\cap S$ invariant, which restricted to $U\cap S$ coincides with $\sD_{f,g}$. Let now $(\hat{U},\hat{z})\in\gU$ be another local chart such that $U\cap\hat{U}\cap S\neq\emptyset$, $\hat{X}_{f,g}$ the corresponding canonical local generator of $\sD_{f,g}$ on $\hat{U}\cap S$ defined as in (\ref{locgen_tang}) and $\hat{\cX}_{f,g}$ its canonical local extension as in (\ref{can_gen_tang_ext}). Being both the charts adapted to $S$ there exists a germ $a\in\cO^*_M$ such that $\hat{z}^1=az^1$. Then one can easily check that 
\begin{equation}\label{can_gen_tang_diff}
\hat{X}_{f,g}=\left(\frac{1}{a|_S}\right)^{\nu} X_{f,g}
\end{equation}
where they overlap. Instead the relation between their extensions is a bit more complicated and we describe it in the next proposition. First, we introduce the following notation:
\begin{enumerate}
 \item[(1)] $T_k$ for holomorphic vector fields of the kind $\sum_{j=2}^n b^j\frac{\partial}{\partial z^j}$, with $b^j\in\cI_S^k$ ({\it tangential terms vanishing on $S$ with `order $k$'});
 \item[(2)] $V_k$ for holomorphic vector fields of the kind $\sum_{j=1}^n c^j\frac{\partial}{\partial z^j}$, with $c^j\in\cI_S^k$ ({\it generic terms vanishing on $S$ with `order $k$'});
\end{enumerate}
\begin{proposition}\label{tang_ext_diff}
Let $M$ be a $n$-dimensional complex manifold, $S\subset M$ a non-singular connected complex hypersurface, $(f,g)\in\text{End}_S^2(M)$ such that $g$ is a local biholomorphism around $S$  and set $\nu=\nu_{f,g}$.\par
Suppose $(f,g)$ tangential along $S$ and let $\gU$ be an atlas adapted to $(g,S)$. If $(U,z)$ and $(\hat{U},\hat{z})$ are two local charts in $\gU$ such that $U\cap\hat{U}\cap S\neq\emptyset$ and $\cX_{f,g}$ and $\hat{\cX}_{f,g}$ are the corresponding local holomorphic vector fields defined as in (\ref{can_gen_tang_ext}), then
$$
\hat{\cX}_{f,g}=\left(\frac{1}{a}\right)^{\nu} \cX_{f,g}+T_{\nu}+V_2
$$
where they overlap, with $a\in\cO^*_M$ such that $\hat{z}^1=az^1$.
\end{proposition}
\begin{pf}
Since $\hat{z}^1=az^1$ then
\begin{equation}\label{frac_1}
\frac{\partial}{\partial\hat{z}^1}=\left(\frac{1}{a}+\frac{\partial\frac{1}{a}}{\partial\hat{z}^1}\hat{z}^1\right)\frac{\partial}{\partial z^1}+\sum_{p=2}^n\frac{\partial z^p}{\partial\hat{z}^1}\frac{\partial}{\partial z^p},
\end{equation}
while
\begin{equation}\label{frac_j}
\frac{\partial}{\partial\hat{z}^j}=\frac{\partial\frac{1}{a}}{\partial\hat{z}^j}\hat{z}^1\frac{\partial}{\partial z^1}+\sum_{p=2}^n\frac{\partial z^p}{\partial\hat{z}^j}\frac{\partial}{\partial z^p},\qquad j=2,\dots,n.
\end{equation}
By (\ref{germs_h}) we have that
\begin{align*}
\hat{h}^1 a^{\nu}(z^1)^{\nu} & =\hat{h}^1(\hat{z}^1)^{\nu}=\hat{w}^1\circ f-\hat{w}^1\circ g= \\
 & =(g_* a\circ f)(w^1\circ f)-(g_* a\circ g)(w^1\circ g)= \\
 & =(w^1\circ f)\left(g_* a\circ f-g_* a\circ g\right)+(g_* a\circ g)h^1(z^1)^{\nu},
\end{align*}
where $w$, $\hat{w}$ are the special coordinates induced by respectively $z$, $\hat{z}$ and we set $g_* a=a\circ g^{-1}$. Then by (\ref{fundamentals}), (\ref{germs_h}) and observing that $w^1\circ f\in\cI_S$ 
$$
\hat{h}^1 a^{\nu}(z^1)^{\nu}=(w^1\circ f)(z^1)^{\nu}\sum_{j=1}^n h^j\frac{\partial a}{\partial z^j}+a h^1(z^1)^{\nu}\phantom{a}\left(\text{mod }\cI_S^{2\nu+1}\right)
$$
since $\frac{\partial g_*a}{\partial w^j}\circ g=\frac{\partial a}{\partial z^j}$. Being $(f,g)$ tangential along $S$ it follows that
$$
w^1\circ f=w^1\circ g \left(\text{mod }\cI_S^{\nu+1}\right)=z^1 \left(\text{mod }\cI_S^{\nu+1}\right),
$$
then we can improve (\ref{htilde_1}) in
\begin{equation}\label{htilde_1_imp}
\hat{h}^1 a^{\nu}=z^1 \sum_{j=2}^n h^j\frac{\partial a}{\partial z^j}+a h^1\phantom{a}\left(\text{mod }\cI_S^2\right).
\end{equation}
Finally, observe that 
\begin{equation}\label{deltcron}
\d_{ij}=\frac{\partial z^i}{\partial z^j}=\sum_{k=1}^n\frac{\partial z^i}{\partial\hat{z}^k}\frac{\partial\hat{z}^k}{\partial z^j},\qquad i,j=1,\dots,n.
\end{equation}
Then manipulating $\hat{\cX}_{f,g}=\sum_{j=1}^n\hat{h}^j\frac{\partial}{\partial\hat{z}^j}$ with (\ref{frac_1}), (\ref{frac_j}), (\ref{htilde_j}), (\ref{htilde_1_imp}) and (\ref{deltcron}) one gets the relation.
\end{pf}
Now assume that $S$ splits into $M$ and fix an atlas $\gU$ on $M$ which is a $\s$-splitting atlas adapted to $(g,S)$. Recall that if $(U,z)\in\gU$ is a local chart such that $U\cap S\neq\emptyset$ then (\ref{locgen_tang_split}) gives the canonical local generator $X_{f,g}^{\s}$ of the foliation $\sD_{f,g}^{\s}$ (eq. $\X_{f,g}^{\s}$) and,  when $\nu=1$, (\ref{locgen_tang_split_nu1}) gives the canonical local generator $X_{f,g}^{\s,1}$ of the foliation $\sD_{f,g}^{\s,1}$ (eq. $\X_{f,g}^{\s,1}$), both associated to these coordinates. Let define $h^1_0=h^1(0,z^2,\dots,z^n)$ (seen as a function on $U$), where $h^1$ is the one appearing in (\ref{germs_h}). Since $z^1$ is a local generator for $\cI_S$ there exists a germ $k^1\in\cO_M$ such that
\begin{equation}\label{hk}
h^1-h^1_0=k^1z^1.
\end{equation}
Then we define the local holomorphic vector fields
\begin{equation}\label{can_gen_comf_emb_ext}
\cX_{f,g}^{\s}=k^1z^1\frac{\partial}{\partial z^1}+\sum_{j=2}^n h^j \frac{\partial}{\partial z^j}
\end{equation}
and, only when $\nu=1$,
\begin{equation}\label{can_gen_comf_emb_ext_nu1}
\cX_{f,g}^{\s,1}=k^1z^1\frac{\partial}{\partial z^1}+(1+h^1_0)\sum_{j=2}^n h^j \frac{\partial}{\partial z^j}
\end{equation}
to be the {\it canonical local extensions} respectively of $X_{f,g}^{\s}$ and $X_{f,g}^{\s,1}$. Observe that both generate a $1$-dimensional holomorphic foliation on $U$ leaving $U\cap S$ invariant which restricted to $U\cap S$ coincides respectively with $\sD_{f,g}^{\s}$ and $\sD_{f,g}^{\s,1}$. Let now $(\hat{U},\hat{z})\in\gU$ be another local chart such that $U\cap\hat{U}\cap S\neq\emptyset$ and $\hat{X}_{f,g}^{\s}$ and $\hat{X}_{f,g}^{\s,1}$ the corresponding canonical local generators of the foliations on $\hat{U}\cap S$ defined as in (\ref{locgen_tang_split}) and (\ref{locgen_tang_split_nu1}). Moreover let $\hat{\cX}_{f,g}^{\s}$ and $\hat{\cX}_{f,g}^{\s,1}$ be their canonical local extensions as defined respectively in (\ref{can_gen_comf_emb_ext}) and (\ref{can_gen_comf_emb_ext_nu1}). As above, there exists a germ $a\in\cO^*_M$ such that $\hat{z}^1=az^1$ and one can easily check that 
\begin{equation}\label{can_gen_comf_emb_diff} 
\hat{X}_{f,g}^{\s}=\left(\frac{1}{a|_S}\right)^{\nu} X_{f,g}^{\s}\quad\text{ and }\quad\hat{X}_{f,g}^{\s,1}=\frac{1}{a|_S}\phantom{a} X_{f,g}^{\s,1}
\end{equation}
where they overlap. We can obtain relations among their local extensions similar to the ones of Proposition \ref{tang_ext_diff} but it is not sufficient that $S$ splits into $M$. In the following we use the same notation of above.
\begin{proposition}\label{comf_emb_ext_diff}
Let $M$ be a $n$-dimensional complex manifold, $S\subset M$ a non-singular connected complex hypersurface, $(f,g)\in\text{End}_S^2(M)$ such that $g$ is a local biholomorphism around $S$  and set $\nu=\nu_{f,g}$.\par
Suppose that $S$ is comfortably embedded into $M$ and let $\gU$ be a $\s$-comfortably atlas adapted to $(g,S)$ (for some splitting morphism $\s$). If $(U,z)$ and $(\hat{U},\hat{z})$ are two local charts in $\gU$ such that $U\cap\hat{U}\cap S\neq\emptyset$, $\cX_{f,g}^{\s}$ and $\hat{\cX}_{f,g}^{\s}$ are the corresponding local holomorphic vector fields defined as in (\ref{can_gen_comf_emb_ext}) and $\cX_{f,g}^{\s,1}$ and $\hat{\cX}_{f,g}^{\s,1}$ the ones defined as in (\ref{can_gen_comf_emb_ext_nu1}), then for $\nu>1$
$$
\hat{\cX}_{f,g}^{\s}=\left(\frac{1}{a}\right)^{\nu} \cX_{f,g}^{\s}+T_1+V_2
$$
where they overlap, with $a\in\cO^*_M$ such that $\hat{z}^1=az^1$. Similarly, if $\nu=1$
$$
\hat{\cX}_{f,g}^{\s,1}=\frac{1}{a}\phantom{a} \cX_{f,g}^{\s,1}+T_1+V_2
$$
where they overlap.
\end{proposition}
\begin{pf}
Let assume $\nu>1$ and focus on the first relation. Recall that since $\gU$ is comfortably then
\begin{equation}\label{first_deriv_null}
\left.\frac{\partial z^j}{\partial\hat{z}^1}\right|_S\equiv \left.\frac{\partial \hat{z}^j}{\partial z^1}\right|_S\equiv 0,\qquad j=2,\dots,n,
\end{equation}
and $\frac{\partial^2 z^1}{(\partial\hat{z}^1)^2}|_S\equiv \frac{\partial^2 \hat{z}^1}{(\partial z^1)^2}|_S\equiv 0$, which implies 
\begin{equation}\label{deriv_a_null}
\left.\frac{\partial a}{\partial z^1}\right|_S \equiv \left.\frac{\partial \frac{1}{a}}{\partial \hat{z}^1}\right|_S \equiv 0.
\end{equation}
Let now $h^1_0$, $k^1$ and $\hat{h}^1_0$, $\hat{k}^1$ be the functions appearing in (\ref{hk}), respectively for the coordinates $z$ and $\hat{z}$. By (\ref{htilde_1}) it follows that $\hat{h}^1_0 a^{\nu}=a h^1_0\phantom{a}(\text{mod }\cI_S)$, since $h^1=h^1_0\phantom{a}(\text{mod }\cI_S)$ and $\hat{h}^1=\hat{h}^1_0\phantom{a}(\text{mod }\cI_S)$. Moreover (\ref{first_deriv_null}) and (\ref{deriv_a_null}) imply that $\frac{\partial (\hat{h}^1_0 a^{\nu}-a h^1_0)}{\partial z^1}|_S\equiv 0$ hence in fact
\begin{equation}\label{interna}
\hat{h}^1_0 a^{\nu}=a h^1_0\phantom{a}\left(\text{mod }\cI_S^2\right).
\end{equation}
Using (\ref{htilde_j}), (\ref{deriv_a_null}) and (\ref{interna}) one has that
\begin{align*}
\hat{k}^1 a^{\nu+1}z^1 & =a^{\nu}\left(\hat{h}^1-\hat{h}^1_0\right)=a^{\nu}\hat{h}^1-a h^1_0\phantom{a}\left(\text{mod }\cI_S^2\right)= \\
& =a h^1-a h^1_0+\sum_{j=2}^n h^j \frac{\partial\hat{z}^1}{\partial z^j} \phantom{a}\left(\text{mod }\cI_S^{\text{min}(\nu,2)}\right)= \\
& =ak^1z^1+z^1\sum_{j=2}^n h^j \frac{\partial a}{\partial z^j}\phantom{a}\left(\text{mod }\cI_S^{\text{min}(\nu,2)}\right),
\end{align*}
then since $\nu>1$ one gets the relation
\begin{equation}\label{k}
\hat{k}^1 a^{\nu+1}=a k^1+\sum_{j=2}^n h^j\frac{\partial a}{\partial z^j}\phantom{a}\left(\text{mod }\cI_S\right),
\end{equation}
which in practice plays here the role of (\ref{htilde_1_imp}). Then manipulating 
$$
\hat{\cX}_{f,g}^{\s}=\hat{k}^1\hat{z}^1\frac{\partial}{\partial\hat{z}^1}+\sum_{j=2}^n\hat{h}^j\frac{\partial}{\partial\hat{z}^j}
$$ 
with the previous (\ref{frac_1}), (\ref{frac_j}), (\ref{htilde_j}), (\ref{deltcron}) and with (\ref{first_deriv_null}), (\ref{deriv_a_null}), (\ref{k}), and recalling that $\nu>1$, one gets the first relation.\par
\smallskip
Let now assume $\nu=1$ and focus on the second relation. Observe that (\ref{k}) is no longer true in this case but we need something similar. Since
\begin{align*}
a\hat{h}^1z^1 & =\hat{h}^1\hat{z}^1=\hat{f}^1-\hat{g}^1= \\
 & =z^1\sum_{j=1}^n h^j \frac{\partial \hat{z}^1}{\partial z^j}+\frac{1}{2}(z^1)^2\sum_{j,k=1}^nh^jh^k \frac{\partial^2\hat{z}^1}{\partial z^j\partial z^k}\phantom{a}\left(\text{mod }\cI_S^3\right)
\end{align*}
then one can easily gets
\begin{align*}
a\hat{h}^1 = ah^1+\left(1+h^1_0\right)z^1\sum_{j=2}^n h^j \frac{\partial a}{\partial z^j} \phantom{a}\left(\text{mod }\cI_S^2\right).
\end{align*}
From this and using (\ref{interna}) again (with $\nu=1$), one obtains
\begin{equation}\label{k_nu1}
a^2 \hat{k}^1=ak^1+\left(1+h^1_0\right)\sum_{j=2}^n h^j\frac{\partial a}{\partial z^j}\phantom{a}\left(\text{mod }\cI_S\right),
\end{equation}
then manipulating 
$$
\hat{\cX}_{f,g}^{\s,1}=\hat{k}^1\hat{z}^1\frac{\partial}{\partial\hat{z}^1}+\left(1+\hat{h}^1_0\right)\sum_{j=2}^n\hat{h}^j\frac{\partial}{\partial\hat{z}^j}
$$ 
with the previous (\ref{frac_1}), (\ref{frac_j}), (\ref{htilde_j}), (\ref{deltcron}), (\ref{first_deriv_null}), (\ref{deriv_a_null}) and with (\ref{k_nu1}) one concludes.
\end{pf}
\begin{remark}
If we consider the couple $(f,\text{Id}_M)$ then Proposition \ref{tang_ext_diff} turns out to be \cite[Lem.4.1.]{ABT2004}, while Proposition \ref{comf_emb_ext_diff} turns out to be \cite[Prop.4.2.]{ABT2004}.
\end{remark}
\bigskip
\section{A Camacho-Sad-type index theorem}\label{CamachoSad}\setcounter{equation}{0}
Let $M$ be a $n$-dimensional complex manifold and $S\subset M$ a globally irreducible complex hypersurface whose regular part is $S'=S-\text{Sing}(S)$. Let $(f,g)\in\text{End}^2_S(M)$ be a couple such that $g$ is a local biholomorphism around $S'$ and with order of coincidence $\nu=\nu_{f,g}$.\par
Let use the following notation in the sequel. When $(f,g)$ is tangential along $S$ set $\sD=\sD_{f,g}$, $\X=\X_{f,g}$ and let $\gU$ denote an atlas adapted to $(g,S')$. If $(U,z)\in\gU$ is a local chart then set also $X=X_{f,g}$ and $\cX=\cX_{f,g}$ as defined in (\ref{locgen_tang}) and (\ref{can_gen_tang_ext}). Instead, when $S'$ is comfortably embedded into $M$ set $\sD=\sD_{f,g}^{\s}$ and $\X=\X_{f,g}^{\s}$ if $\nu>1$ while set $\sD=\sD_{f,g}^{\s,1}$ and $\X=\X_{f,g}^{\s,1}$ if $\nu=1$. Moreover let $\gU$ denote a $\s$-comfortably atlas adapted to $(g,S')$ (for some splitting morphism $\s$) and if $(U,z)\in\gU$ is a local chart then set $X=X_{f,g}^{\s}$ and $\cX=\cX_{f,g}^{\s}$ as defined in (\ref{locgen_tang_split}) and (\ref{can_gen_comf_emb_ext}) when $\nu>1$, while set $X=X_{f,g}^{\s,1}$ and $\cX=\cX_{f,g}^{\s,1}$ as defined in (\ref{locgen_tang_split_nu1}) and (\ref{can_gen_comf_emb_ext_nu1}) when $\nu=1$. In all the cases let denote $S^0=S'-\text{Sing}(\sD)$ and $\sD^0=\sD|_{S^0}$.\par
Let $(U,z)\in\gU$ be such that $U\cap S^0\neq\emptyset$. As already observed the local holomorphic vector field $\cX$ generates a $1$-dimensional holomorphic foliation on $U$ which extends $\sD|_{U\cap S'}$, then locally we are in the situation described at the beginning of Section \ref{loc_ext}. This means that for any $(U,z)\in\gU$ (such that $U\cap S^0\neq\emptyset$) we can define a partial holomorphic connection on $N_{U\cap S^0}$ along $N_{U\cap S^0}^{\ot\nu}\equiv \X|_{U\cap S^0}\subset TS'|_{U\cap S^0}$ as in (\ref{CS_conn}), that is
\begin{align}\label{CS_conn_loc}
\cN_{U\cap S^0} & \stackrel{\d_U^{cs}}{\longrightarrow} \left(\cN_{U\cap S^0}^{\ot\nu}\right)^*\ot\cN_{U\cap S^0} \notag \\
w & \longrightarrow \d_U^{cs}(w)\text{ s.t. }\d_U^{cs}(w)(\f X)=\f\pi\left(\left.[\cX,\tilde{w}]\right|_{S^{0}}\right)
\end{align}
for any $w\in\cN_{U\cap S^0}$ and $\f\in\cO_{U\cap S^0}$, where $\pi:TM|_{S^0}\to N_{S^0}$ is the projection and $\tilde{w}\in\cT_M|_{U\cap S^0}$ is any vector field such that $\pi(\tilde{w}|_{U\cap S^{0}})=w$. Observe that since $X$ is a local generator of $\sD$ on $U\cap S'$ then any local vector field $v\in\cN_{U\cap S^0}^{\ot\nu}$ is of the form $v=\f X$ for some function $\f\in\cO_{U\cap S^0}$.\par
The key point is that we can glue together all these local partial holomorphic connections by Proposition \ref{tang_ext_diff} and \ref{comf_emb_ext_diff} and consequently define a (global) partial holomorphic connection on $N_{S^0}$ along $N_{S^0}^{\ot\nu}\equiv \X|_{S^0}\subset TS^0$.
\begin{prop}\label{CS_conn_new}
Let $M$ be a $n$-dimensional complex manifold, $S\subset M$ a globally irreducible complex hypersurface and $(f,g)\in\text{End}^2_S(M)$ a couple such that $g$ is a local biholomorphism around $S'$ with order of coincidence $\nu=\nu_{f,g}$. Assume that $(f,g)$ is tangential along $S$ or that $S'$ is comfortably embedded into $M$ and use the notation introduced above.\par
Then if $(U,z)$ and $(\hat{U},\hat{z})$ are two local charts in $\gU$ such that $U\cap\hat{U}\cap S^0\neq\emptyset$ and $\d_U^{cs}$ and $\d_{\hat{U}}^{cs}$ are the corresponding local partial holomorphic connections defined as in (\ref{CS_conn_loc}) we have that
$$
\d_U^{cs}=\d_{\hat{U}}^{cs}
$$
where they overlap. Consequently, there is a well-defined partial holomorphic connection 
$$
\cN_{S^0} \stackrel{\d^{cs}_{\sD}}{\longrightarrow} \left(\cN_{S^0}^{\ot\nu}\right)^*\ot\cN_{S^0}
$$
on $N_{S^0}$ along $N_{S^0}^{\ot\nu}$.
\end{prop}
\begin{pf}
Since $X$ is a local generator of $\sD$ on $U\cap\hat{U}\cap S^0$ we just need to prove that
$$
\d_{\hat{U}}^{cs}(w)(X)=\d_{U}^{cs}(w)(X)
$$
for any $w\in\cN_{U\cap\hat{U}\cap S^0}$. This is true because if $a\in\cO^*_M$ is such that $\hat{z}^1=az^1$ and $\pi:TM|_{S^0}\to N_{S^0}$ is the projection then by  (\ref{can_gen_tang_diff}), (\ref{can_gen_comf_emb_diff}) and by Proposition \ref{tang_ext_diff} and \ref{comf_emb_ext_diff} it follows that
\begin{align*}
\d_{\hat{U}}^{cs}(w)(X) & =\left(\left.a\right|_{S^0}\right)^{\nu}\d_{\hat{U}}^{cs}(w)(\hat{X})=\left(\left.a\right|_{S^0}\right)^{\nu}\pi\left(\left.\left[\hat{\cX},\tilde{w}\right]\right|_{S^0}\right)= \\
 & =\left(\left.a\right|_{S^0}\right)^{\nu}\pi\left(\left.\left[\left(\frac{1}{a}\right)^{\nu}\cX+T_1+V_2,\tilde{w}\right]\right|_{S^0}\right)= \\
 & =\pi\left(\left.\left[\cX,\tilde{w}\right]\right|_{S^0}\right)=\d_{U}^{cs}(w)(X),
\end{align*}
where the last but one equality is due to the fact that $[T_1,\tilde{w}]|_{S^0}\in\cT_{S^0}$ and $[V_2,\tilde{w}]|_{S^0}\equiv 0$.
\end{pf}
Now recall that a hypersurface $S\subset M$ defines naturally a line bundle on $M$, usually denoted by $[S]$, and that $[S]|_{S'}\cong N_{S'}$. Moreover $S$ can be seen as the zero locus of a holomorphic section of $[S]$. Then by Proposition \ref{CS_conn_new} and by \cite[Th.IV.2.4.]{Su1998} (or \cite[Th.5.3.7.]{BSS2009}) we have the following theorem.
\begin{theorem}[Camacho-Sad-type index theorem]\label{CS_theorem}
Let $M$ be a $n$-dimensional complex manifold, $S\subset M$ a globally irreducible compact complex hypersurface whose regular part is $S'=S-\text{Sing}(S)$ and $(f,g)\in\text{End}^2_S(M)$ a couple such that $g$ is a local biholomorphism around $S'$ with order of coincidence $\nu=\nu_{f,g}$. Assume that
\begin{itemize}
\item[(i)] $(f,g)$ is tangential along $S$ 
\end{itemize}
or that
\begin{itemize}
\item[(ii)] $S'$ is comfortably embedded into $M$.
\end{itemize}
In case $(i)$ set $\sD=\sD_{f,g}$ while in case $(ii)$ set $\sD=\sD_{f,g}^{\s}$ when $\nu>1$ or $\sD=\sD_{f,g}^{\s,1}$ when $\nu=1$ (where $\s$ is some splitting morphism) and suppose $\sD\neq 0$. Let $\text{Sing}(S)\cup\text{Sing}(\sD)=\sqcup_{\l}\S_{\l}$ be the decomposition in connected components of the singular set $\text{Sing}(S)\cup\text{Sing}(\sD)$. \par
Then there exist complex numbers $\text{Res}(\sD;S;\S_{\l})$ such that
$$
\sum_{\l}\text{Res}\left(\sD;S;\S_{\l}\right)=\int_S c_1^{n-1}\left([S]\right),
$$
where $c_1([S])$ denotes the first Chern class of the line bundle $[S]$.
\end{theorem}
We call Theorem \ref{CS_theorem} a {\it Camacho-Sad-type index theorem} because it is basically inspired by the Camacho-Sad index theorem proved in \cite[Appendix]{CS1982}. 
To be precise if $n=2$, $S$ is a non-singular curve and $g=\text{Id}_M$ then Theorem \ref{CS_theorem} is broadly the residual index theorem \cite[Th.1.1.]{ab2001}, which was inspired by the Camacho-Sad index theorem just mentioned. 
See the survey \cite{ab2008} for a clear exposition, without too many details, of a general procedure to obtain Camacho-Sad-type index theorems. A full explanation of this procedure and some resulting theorems may be found in \cite{ABT2008}.
\begin{remark}
If we consider the couple $(f,\text{Id}_M)$ then Theorem \ref{CS_theorem} turns out to be \cite[Th.6.2.]{ABT2004}, which itself generalize \cite[Th.1.1.]{ab2001}. Observe that in their notation $\text{Sing}(f,\text{Id}_M)$ is $\text{Sing}(f)$ and they assume $S'-\text{Sing}(f)$  to be comfortably embedded in $M$, not $S'$. In fact their assumption is sufficient. 
\end{remark}
We conclude this section by deriving explicit formulas for the computation of the residues in Theorem \ref{CS_theorem} at isolated singular points. First, we briefly recall how these residues are defined in Lehmann-Suwa theory. As a reference see \cite[Sec.IV.2. and IV.6.]{Su1998}.\par
Let $U^0\subset M$ be a tubular neighborhood of $S^0$ (it always exists, see for example \cite[p.465]{Sp1979}) and $r:U^0\to S^0$ the associated smooth retraction. Observe that $r^*N_{S^0}\cong [S]|_{U^0}$ hence if $\n_{\sD}^{cs}$ is any $(1,0)$-type extension of $\d_{\sD}^{cs}$ then $r^*\n_{\sD}^{cs}$ is a smooth connection on $[S]|_{U^0}$ such that, by the Bott vanishing theorem,
$$
c_1^{n-1}\left(r^*\n_{\sD}^{cs}\right)=r^*c_1^{n-1}\left(\n_{\sD}^{cs}\right)=0.
$$
Let $\S$ be a connected component of $\text{Sing}(S)\cup\text{Sing}(\sD)$, $U\subset M$ any open neighborhood of it such that $U\cap(\text{Sing}(S)\cup\text{Sing}(\sD))=\S$ and take any connection $\n_{U}$ on $[S]|_U$. Lastly, let $\tilde{R}\subset U$ be any  compact real sub-manifold of dimension $2n$ oriented as $M$ and such that $\S\subset\text{int}(\tilde{R})$. Assume its boundary $\partial\tilde{R}$ transverse to $S$ and with the orientation induced by $\tilde{R}$. Set $R=\tilde{R}\cap S$ and $\partial R=\partial\tilde{R}\cap S$. Then by definition the residue is 
$$
\text{Res}\left(\sD;S;\S\right)=\int_R c_1^{n-1}\left(\n_U\right)-\int_{\partial R}c_1^{n-1}\left(r^*\n_{\sD}^{cs},\n_U\right),
$$
where $c_1^{n-1}(r^*\n_{\sD}^{cs},\n_U)$ is the Bott difference form. One can show that this formula does not depend on the various choices done.\par
In particular, if $\S=\{p\}$ we can suppose $U$ to be such that $[S]|_U$ is trivial and then take $\n_{U}$ trivial respect to some generator. With this choice
$$
\text{Res}\left(\sD;S;p\right)=-\int_{\partial R}c_1^{n-1}\left(r^*\n_{\sD}^{cs},\n_U\right)
$$
and since $\partial R\subset U\cap S^0$ and $(r^*\n_{\sD}^{cs})|_{S^0}=\n_{\sD}^{cs}$ we can in fact rewrite it as
$$
\text{Res}\left(\sD;S;p\right)=-\int_{\partial R}c_1^{n-1}\left(\n_{\sD}^{cs},\left.\n_U\right|_{U\cap S^0}\right).
$$
Observe that $\n_U|_{U\cap S^0}$ is a smooth connection on $N_{U\cap S^0}$.\par
\smallskip
Assume now that $p\in\text{Sing}(\sD)-\text{Sing}(S)$. In this case $N_{U\cap S}$ does exist and we can substitute $\n_U|_{U\cap S^0}$ with a smooth connection $\n_{U\cap S}$ on $N_{U\cap S}$ trivial respect to some local generator of it, thus 
$$
\text{Res}\left(\sD;S;p\right)=-\int_{\partial R}c_1^{n-1}\left(\n_{\sD}^{cs},\n_{U\cap S}\right).
$$
At this point one can argue as in the proof of \cite[Th.III.5.5.]{Su1998} and obtain explicit formulas. In particular, let $(U,z)$ be a local chart of $M$ at $p$ belonging to $\gU$ and the $h^j\in\cO_{M,p}$ as in (\ref{germs_h}). Choose as $\n_{U\cap S}$ the connection trivial respect to $\partial z^1$. Then when $(f,g)$ is tangential along $S$ and $\sD=\sD_{f,g}$ we get the formula 
\begin{equation}\label{CS_formula_1}
\text{Res}\left(\sD;S;p\right)=\left(\frac{-i}{2\pi}\right)^{n-1}\int_{\G}\frac{\left(\ell^1\right)^{n-1}}{h^2\cdots h^n}\text{ d}z^2\wedge\cdots\wedge\text{d}z^n,
\end{equation}
where $\ell^1\in\cO_{M,p}$ is the one such that $h^1=\ell^1z^1$ and $\G$ is the $(n-1)$ cycle
$$
\G=\left\{q\in U\cap S\text{ s.t. }|h^2(q)|=\dots=|h^n(q)|=\e\right\},
$$
for $\e>0$ small enough, oriented so that $\text{d}\q^2\wedge\cdots\wedge\text{d}\q^n>0$ where $\q^j=\text{arg}(h^j)$ for $j=2,\dots,n$. When $S'$ is comfortably embedded into $M$, $\nu>1$ and $\sD=\sD_{f,g}^{\s}$ we get instead the formula 
\begin{equation}\label{CS_formula_2}
\text{Res}\left(\sD;S;p\right)=\left(\frac{-i}{2\pi}\right)^{n-1}\int_{\G}\frac{\left(k^1\right)^{n-1}}{h^2\cdots h^n}\text{ d}z^2\wedge\cdots\wedge\text{d}z^n,
\end{equation}
where $k^1$ is the one appearing in (\ref{hk}) and $\G$ is as above. Finally, when $S'$ is comfortably embedded into $M$, $\nu=1$ and $\sD=\sD_{f,g}^{\s,1}$ we have 
\begin{equation}\label{CS_formula_3}
\text{Res}\left(\sD;S;p\right)=\left(\frac{-i}{2\pi}\right)^{n-1}\int_{\G'}\frac{\left(k^1\right)^{n-1}}{\left(1+h^1_0\right)^{n-1}h^2\cdots h^n}\text{ d}z^2\wedge\cdots\wedge\text{d}z^n,
\end{equation}
where $h^1_0$ is the one appearing in (\ref{hk}) and 
$$
\G'=\left\{q\in U\cap S\text{ s.t. }\left|(1+h^1(q))h^2(q)\right|=\dots=\left|(1+h^1(q))h^n(q)\right|=\e\right\}.
$$
\begin{remark}\label{calc_res_altern}
Let $p\in\text{Sing}(\sD)-\text{Sing}(S)$ and $(U,z)\in\gU$ be a local chart at $p$. By construction $\d^{cs}_{\sD}$ is locally induced by the natural `holomorphic action' of $X$ on $N_{U\cap S^0}$ by Lie bracket then we could have used directly the formula just after \cite[Th.2]{LS1995} (which is the same of \cite[Th.IV.6.3.]{Su1998}) to get (\ref{CS_formula_1}), (\ref{CS_formula_2}) and (\ref{CS_formula_3}). See also \cite[Rmk.IV.6.7.(1) and (4)]{Su1998}.
\end{remark}
We end treating the case $p\in\text{Sing}(S)$, which is more complicated since $N_{U\cap S}$ does not exist and in general there is no natural local extension of $\sD$ at such a point. Anyway, assuming that $g$ is a biholomorphism in a neighborhood of $p$ we can calculate explicitly the residue when $n=2$ and, in some cases, when $n>2$.\par
Since $p\in\text{Sing}(S)$ we do not have a local chart $(U,z)\in\gU$ at $p$ but when $n=2$ we can `almost' find it. In fact, let $U\subset M$ be an open neighborhood  of $p$ such that $U\cap(\text{Sing}(\sD)\cup\text{Sing}(S))=\{p\}$ and $g|_U$ is a biholomorphism onto its image. We can assume to have on it a local generator $y$ of $\cI_S$ and some coordinates $(u^1,u^2)$ such that $\text{d}y\wedge \text{d}u^2\neq 0$ on $U\cap S^0$. In particular $U\cap S=\{y=0\}$ and we can suppose (possibly shrinking $U$) that $(U-\{p\},(y,u^2))$ is a local chart in $\gU$. Since $y$ generates $\cI_S$ on $U$, the dual of $[y]\in\cI_S / \cI_S^2$ defines a local generator of $[S]$ on $U$ whose restriction to $U\cap S^0$ coincides with the local generator $\partial y=\pi (\frac{\partial}{\partial y}|_{U\cap S^0})$ of $N_{U\cap S^0}$. We can then choose $\n_U$ to be trivial respect to the dual of $[y]$ and consequently $\n_U|_{U\cap S^0}$ is trivial respect to $\partial y$. With this choice we can again argue as in the proof of \cite[Th.III.5.5.]{Su1998} and obtain explicit formulas likewise the previous ones. Referring to the notation of (\ref{CS_formula_1}), (\ref{CS_formula_2}) and (\ref{CS_formula_3}) and reminding how are defined the special coordinates $(w^1,w^2)$ associated to the $(z^1,z^2)$ observe that when $(f,g)$ is tangential along $S$ 
$$
\left.\frac{\ell^1}{h^2}\right|_{S^0}=\left.\frac{f^1-g^1}{z^1(f^2-g^2)}\right|_{S^0}=\left.\frac{z^1\circ g^{-1}\circ f-z^1}{z^1(z^2\circ g^{-1}\circ f-z^2)}\right|_{S^0}.
$$
Instead when $S'$ is comfortably embedded into $M$, if $\nu>1$
$$
\left.\frac{k^1}{h^2}\right|_{S^0}=\left.\frac{f^1-g^1-(z^1)^{\nu}h^1_0}{z^1(f^2-g^2)}\right|_{S^0}=\left.\frac{z^1\circ g^{-1}\circ f-z^1-(z^1)^{\nu}h^1_0}{z^1(z^2\circ g^{-1}\circ f-z^2)}\right|_{S^0},
$$
while if $\nu=1$
\begin{align*}
\left.\frac{k^1}{(1+h^1_0)h^2}\right|_{S^0} & =\left.\frac{f^1-g^1-z^1h^1_0}{(z^1+z^1h^1_0)(f^2-g^2)}\right|_{S^0}=\\
 & = \left.\frac{z^1\circ g^{-1}\circ f-z^1-z^1h^1_0}{(z^1\circ g^{-1}\circ f)(z^2\circ g^{-1}\circ f-z^2)}\right|_{S^0},
\end{align*}
where we use the fact that $h^1_0|_{S^0}=h^1|_{S^0}$. Hence substituting the coordinates $(z^1,z^2)$ with the $(y,u^2)$ and using these equalities we have the following. When $(f,g)$ is tangential along $S$ and $\sD=\sD_{f,g}$ by (\ref{CS_formula_1}) we get
\begin{equation}\label{CS_formula_4}
\text{Res}\left(\sD;S;p\right)=\frac{1}{2\pi i}\int_{\g}\frac{y\circ g^{-1}\circ f-y}{y(u^2\circ g^{-1}\circ f-u^2)}\text{ d}u^2,
\end{equation}
where $\g\subset U\cap S^{\circ}$ is any simple closed curve around $p$ with the orientation induced by the one of $S'$. When $S'$ is comfortably embedded into $M$, $\nu>1$ and $\sD=\sD_{f,g}^{\s}$ by (\ref{CS_formula_2}) we get instead
\begin{equation}\label{CS_formula_5}
\text{Res}\left(\sD;S;p\right)=\frac{1}{2\pi i}\int_{\g}\frac{y\circ g^{-1}\circ f-y-y^{\nu}b}{y(u^2\circ g^{-1}\circ f-u^2)}\text{ d}u^2,
\end{equation}
where $b=\frac{y\circ g^{-1}\circ f-y}{y^{\nu}}$ and $\g$ as before. Finally, if $S'$ is comfortably embedded into $M$, $\nu=1$ and $\sD=\sD_{f,g}^{\s,1}$ by (\ref{CS_formula_3}) we obtain
\begin{equation}\label{CS_formula_6}
\text{Res}\left(\sD;S;p\right)=\frac{1}{2\pi i}\int_{\g}\frac{y\circ g^{-1}\circ f-y-y b}{(y\circ g^{-1}\circ f)(u^2\circ g^{-1}\circ f-u^2)}\text{ d}u^2,
\end{equation}
with $b$ and $\g$ as before.\par
If $n>2$ we are able to compute explicitly the residue at $p\in\text{Sing}(S)$ when $(f,g)$ is tangential along $S$ and $\nu>1$. In fact, with these hypotheses we can define a local holomorphic vector field around $p$ leaving $S$ invariant with (possibly) an isolated singularity at $p$, whose natural `holomorphic action' produced by Lie bracket on $N_{S^0}$ (see \cite[Sec.IV.6.]{Su1998}) induces locally $\d^{cs}_{\sD}$. In this way we can apply directly the formula just after \cite[Th.2]{LS1995}.\par
For this purpose, let $U\subset M$ be an open neighborhood of $p$ such that $U\cap(\text{Sing}(S)\cup\text{Sing}(sD))=\{p\}$ and $g|_U$ is a biholomorphism onto its image. Suppose to have a local generator $y$ of $\cI_S$ and any coordinates $(u^1,\dots,u^n)$ on it, with $(v^1=u^1\circ g^{-1},\dots,v^n=u^n\circ g^{-1})$ the corresponding ``special'' coordinates on $g(U)$ at $f(p)$.  Let define the holomorphic vector field on $U$
\begin{equation}\label{vector_V}
\cV_{f,g}=\sum_{j=1}^n\frac{v^j\circ f-v^j\circ g}{y^{\nu}}\frac{\partial}{\partial u^j}
\end{equation}
and $U^j=\{x\in U\text{ s.t. }\frac{\partial y}{\partial u^j}(x)\neq 0\}$ for $j=1,\dots,n$. Observe that $\sqcup_{j=1}^nU^j$ is an open subset of $U-\{p\}$ containing $U\cap S^0$ and in particular the sets $U^j\cap S^0$ cover $U\cap S^0$. We can put on each $U^j$ the coordinates $z_j=(z_j^1,\dots,z_j^n)$ adapted to $(g,S')$ defined by
\begin{align*}
z^1_{j}&=y, \\
z^i_{j}&=u^{i-1},&\text{ for } i=2,\dots,j, \\
z^i_{j}&=u^i,&\text{ for } i=j+1,\dots,n.
\end{align*}
Clearly $(U^j,z_j)\in\gU$ for every $j=1,\dots,n$ and let $w_j=(w_j^1,\dots,w_j^n)$ be the usual special coordinates associated to the $z_j$. If $\cX^j$ is the local holomorphic vector field associated to the chart $(U^j,z_j)$ defined as in (\ref{can_gen_tang_ext}) then one can prove that
\begin{equation}\label{vectors_V_X}
\left.\cV_{f,g}\right|_{U^j} = \cX^j+V_{\nu}^j,\qquad\quad\text{for }j=1,\dots,n,
\end{equation}
where $V_{\nu}^j$ is a holomorphic vector field on $U^j$ whose coefficients are in $\cI_S^{\nu}$. We show (\ref{vectors_V_X}) only for $j=1$ since the other cases can be proved in the same way. For simplicity let denote the coordinates by $z=(z^1,\dots,z^n)$ and $w=(w^1,\dots,w^n)$. Let $h^1,\dots,h^n$ be the corresponding germs as in (\ref{germs_h}), then on $U^1$
$$
\frac{v^i\circ f-v^i\circ g}{y^{\nu}}=\frac{w^i\circ f-w^i\circ g}{(z^1)^{\nu}}=h^i
$$
for $i=2,\dots,n$, while reminding (\ref{fundamentals}) 
\begin{align*}
v^1\circ f-v^1\circ g=\left(\frac{\partial y}{\partial u^1}\right)^{-1}h^1 y^{\nu} -\left(\frac{\partial y}{\partial u^1}\right)^{-1}\sum_{k=2}^n\frac{\partial y}{\partial u^k}h^k y^{\nu}\left(\text{mod }\cI^{2\nu}_S\right)
\end{align*}
hence
$$
\frac{v^1\circ f-v^1\circ g}{y^{\nu}}=\left(\frac{\partial y}{\partial u^1}\right)^{-1}\left[h^1-\sum_{k=2}^n\frac{\partial y}{\partial u^k}h^k\right]\left(\text{mod }\cI^{\nu}_S\right).
$$
Therefore putting these equalities into (\ref{vector_V}) and noting that $\frac{\partial}{\partial u^1}=\frac{\partial y}{\partial u^1}\frac{\partial}{\partial z^1}$ and $\frac{\partial}{\partial u^i}=\frac{\partial y}{\partial u^i}\frac{\partial}{\partial z^1}+\frac{\partial}{\partial z^i}$ for $i=2,\dots,n$ we done.\par
By (\ref{vectors_V_X}) it follows that $\cV_{f,g}$ generates a $1$-dimensional holomorphic foliation on $U$ leaving $U\cap S^0$ invariant, which restricted to $U\cap S^0$ coincides with $\sD_{f,g}$ and having an isolated singularity at $p$. Moreover taking into account that $\nu>1$ and reminding (\ref{CS_conn_loc}), equation (\ref{vectors_V_X}) shows even that the natural `holomorphic action' of $\cV_{f,g}$ on $N_{U\cap S^0}$ induces locally $\d^{cs}_{\sD}$. Hence choosing the coordinates $(u^1,\dots,u^n)$ in such a way that $\{y=(v^2\circ f-v^2\circ g)/y^{\nu}=\cdots=(v^n\circ f-v^n\circ g)/y^{\nu}=0\}=\{p\}$ (we can do that, see for example \cite[Cor.IV.4.5.]{Su1998}) then we can express the residue by the formula just after \cite[Th.2]{LS1995} (eq. \cite[Th.IV.6.3.]{Su1998}). In particular, observing that
$$
\cV_{f,g}\left(y\right)=\left(\sum_{j=1}^n\frac{v^j\circ f-v^j\circ g}{y^{\nu+1}}\frac{\partial y}{\partial u^j} \right)y
$$
we get the formula
\begin{equation}\label{CS_formula_7}
\text{Res}\left(\sD;S;p\right)=\left(\frac{-i}{2\pi}\right)^{n-1}\int_{\G} \frac{\left[\sum_{j=1}^n(v^j\circ f-v^j\circ g)\frac{\partial y}{\partial u^j}\right]^{n-1}}{y^{n-1}\prod_{j=2}^n(v^j\circ f-v^j\circ g)} \text{ d}u^2\wedge\cdots\wedge\text{d}u^n,
\end{equation}
where $\G$ is the $(n-1)$ cycle
$$
\G=\left\{q\in U\cap S\text{ s.t. }\left|\frac{v^j\circ f-v^j\circ g}{y^{\nu}}(q)\right|=\e, \text{ for }j=2,\dots,n\right\},
$$
for $\e>0$ small enough, oriented as usual.\par
\begin{remark}
If $n=2$, $(f,g)$ is tangential along $S$ and $\nu>1$ we can argue as just done, taking as local chart $(U,u)$ at $p$ one such that $\text{d}y\wedge\text{d}u^2\neq 0$ on $U\cap S^0$. Thus we can assume (possibly shrinking $U$) that $U-\{p\}=U^1$ and then that $(U-\{p\},(y,u^2))\in\gU$. If $\cX$ is the local holomorphic vector field associated to $(y,u^2)$ as in (\ref{can_gen_tang_ext}) (defined on $U-\{p\}$) then by (\ref{vectors_V_X}) the natural `holomorphic actions' of $\cV_{f,g}$ and $\cX$ on $N_{U\cap S^0}$ are in fact the same. Then suitably modifying (\ref{CS_formula_7}) we recover (\ref{CS_formula_4}).\par
Lastly, observe that if $p$ was a non-singular point for $S$ then we could take a local chart $(U,u)=(U,z)\in\gU$ and as local generator of $\cI_S$ the function $y=z^1$. With these choices clearly $\cV_{f,g}=\cX$ and then (\ref{CS_formula_7}) would be (\ref{CS_formula_1}).
\end{remark}
\bigskip
\section{A Lehmann-Suwa-type index theorem}\label{LehmannSuwa}\setcounter{equation}{0}
Let $M$, $S$ and $(f,g)$ be as in Section \ref{CamachoSad} and set as usual $\nu=\nu_{f,g}$ and $S'=S-\text{Sing}(S)$. From now on assume $(f,g)$ tangential along $S$ and let denote $\sD=\sD_{f,g}$,  $S^0=S'-\text{Sing}(\sD)$ and $\X=\X_{f,g}$ for simplicity. Moreover let $\gU$ be an atlas adapted to $(g,S')$ and if $(U,z)\in\gU$ set also $X=X_{f,g}$ and $\cX=\cX_{f,g}$, as defined respectively in (\ref{locgen_tang}) and (\ref{can_gen_tang_ext}).\par
Let $(U,z)\in\gU$ be such that $U\cap S^0\neq\emptyset$. As in Section \ref{CamachoSad} the local holomorphic vector field $\cX$ generates a $1$-dimensional holomorphic foliation on $U$ which extends $\sD|_{U\cap S'}$, hence we can define a partial holomorphic connection on $N_{\sD}^M|_{U\cap S^0}$ along $N_{U\cap S^0}^{\ot\nu}\equiv \X|_{U\cap S^0}\subset TS'|_{U\cap S^0}$ as in (\ref{LS_conn}), that is
\begin{align}\label{LS_conn_loc}
\cO\left(\left.N_{\sD}^M\right|_{U\cap S^0}\right) & \stackrel{\d_U^{ls}}{\longrightarrow} \left(\cN_{U\cap S^0}^{\ot\nu}\right)^*\ot\cO\left(\left.N_{\sD}^M\right|_{U\cap S^0}\right) \notag \\
w & \longrightarrow \d_U^{ls}(w)\text{ s.t. }\d_U^{ls}(w)(\f X)=\f\r\left(\left.[\cX,\tilde{w}]\right|_{S^{0}}\right)
\end{align}
for any $w\in\cO(N_{\sD}^M|_{U\cap S^0})$ and $\f\in\cO_{U\cap S^0}$, where $\r:TM|_{S^0}\to N_{\sD}^M$ is the projection and $\tilde{w}\in\cT_M|_{U\cap S^0}$ is any vector field such that $\r(\tilde{w}|_{U\cap S^{0}})=w$. Observe that since $X$ is a local generator of $\sD$ on $U\cap S'$ then any local vector field $v\in\cN_{U\cap S^0}^{\ot\nu}$ is of the form $v=\f X$ for some $\f\in\cO_{U\cap S^0}$.\par
By Proposition \ref{tang_ext_diff} we can glue together all these local partial holomorphic connections {\it when $\nu>1$} and then define a (global) partial holomorphic connection on $N_{\sD}^M$ along $N_{S^0}^{\ot\nu}\equiv\X|_{S^0}\subset TS^0$).
\begin{prop}\label{LS_conn_new}
Let $M$ be a $n$-dimensional complex manifold and $S\subset M$ a globally irreducible complex hypersurface. Let $(f,g)\in\text{End}^2_S(M)$ be a couple such that $g$ is a local biholomorphism around $S'$, with order of coincidence $\nu=\nu_{f,g}>1$ and tangential along $S$.\par
Then, referring to the notation introduced above, if $(U,z)$ and $(\hat{U},\hat{z})$ are two local charts in $\gU$ such that $U\cap\hat{U}\cap S^0\neq\emptyset$ and $\d_U^{ls}$ and $\d_{\hat{U}}^{ls}$ are the corresponding local partial holomorphic connections defined as in (\ref{LS_conn_loc}) we have that
$$
\d_U^{ls}=\d_{\hat{U}}^{ls}
$$
where they overlap. Consequently, there is a well-defined partial holomorphic connection 
$$
\cO\left(N_{\sD}^M\right) \stackrel{\d^{ls}_{\sD}}{\longrightarrow} \left(\cN_{S^0}^{\ot\nu}\right)^*\ot\cO\left(N_{\sD}^M\right)
$$
on $N_{\sD}^M$ along $N_{S^0}^{\ot\nu}$.
\end{prop}
\begin{pf}
Since $X$ is a local generator of $\sD$ on $U\cap\hat{U}\cap S^0$ we just need to prove that
$$
\d_{\hat{U}}^{ls}(w)(X)=\d_{U}^{ls}(w)(X)
$$
for any $w\in\cO(N_{\sD}^M|_{U\cap\hat{U}\cap S^0})$. This is true because if $a\in\cO^*_M$ is such that $\hat{z}^1=az^1$ and $\r:TM|_{S^0}\to N_{\sD}^M$ is the projection then by  (\ref{can_gen_tang_diff}) and Proposition \ref{tang_ext_diff}, recalling that $\nu>1$, it follows 
\begin{align*}
\d_{\hat{U}}^{ls}(w)(X) & =\left(\left.a\right|_{S^0}\right)^{\nu}\d_{\hat{U}}^{ls}(w)(\hat{X})=\left(\left.a\right|_{S^0}\right)^{\nu}\r\left(\left.\left[\hat{\cX},\tilde{w}\right]\right|_{S^0}\right)= \\
 & =\left(\left.a\right|_{S^0}\right)^{\nu}\r\left(\left.\left[\left(\frac{1}{a}\right)^{\nu}\cX+V_2,\tilde{w}\right]\right|_{S^0}\right)= \\
 & =\r\left(\left.\left[\cX,\tilde{w}\right]\right|_{S^0}\right)=\d_{U}^{ls}(w)(X),
\end{align*}
where the last but one equality is due to the fact that $[V_2,\tilde{w}]|_{S^0}\equiv 0$.
\end{pf}
Recall again that a hypersurface $S\subset M$ defines naturally a line bundle $[S]$ on $M$ such that $[S]|_{S'}\cong N_{S'}$. Then by Proposition \ref{LS_conn_new} and arguing as at \cite[p.130]{Su1998} we get the following theorem.
\begin{theorem}[Lehmann-Suwa-type index theorem]\label{LS_theorem}
Let $M$ be a $n$-dimensional complex manifold and $S\subset M$ a globally irreducible compact complex hypersurface. Let $(f,g)\in\text{End}^2_S(M)$ be a couple such that $g$ is a local biholomorphism around $S'$, with order of coincidence $\nu=\nu_{f,g}>1$ and tangential along $S$. Set $\sD=\sD_{f,g}$ and let $\text{Sing}(S)\cup\text{Sing}(\sD)=\sqcup_{\l}\S_{\l}$ be the decomposition in connected components of the singular set $\text{Sing}(S)\cup\text{Sing}(\sD)$.\par
Then for any symmetric homogeneous polynomial $\f\in\bC[z_1,\dots,z_{n-1}]$ of degree $n-1$ there exist complex numbers $\text{Res}_{\f}(\sD;TM|_S-[S]|_S^{\ot\nu};\S_{\l})$ such that
$$
\sum_{\l}\text{Res}\left(\sD;\left.TM\right|_S-\left.[S]\right|_S^{\ot\nu};\S_{\l}\right)=\int_S \f\left(TM-[S]^{\ot\nu}\right).
$$
\end{theorem}
We call Theorem \ref{LS_theorem} a {\it Lehmann-Suwa-type index theorem} since Lehmann and Suwa have introduced this kind of residues and the partial connection (\ref{LS_conn}), which inspires the one in Proposition \ref{LS_conn_new} (see \cite{LS1995}, \cite{LS1999} and also \cite{KS1997}). 
\begin{remark}
If we consider the couple $(f,\text{Id}_M)$ then Theorem \ref{LS_theorem} turns out to be \cite[Th.6.3.]{ABT2004}. Observe that they assume also $S'$ comfortably embedded into $M$ but it is not necessary.
\end{remark}
Like in Section \ref{BaumBott} and \ref{CamachoSad}, we conclude by deriving explicit formulas for the residues in Theorem \ref{LS_theorem} at isolated singularities. As usual we briefly recall how the residues are defined in Lehmann-Suwa theory (for a reference see \cite[Sec.IV.5.]{Su1998}).\par
This kind of residues are defined very similarly to the ones in Section \ref{BaumBott}. Let $\f\in\bC[z_1,\dots,z_{n-1}]$ be any symmetric homogeneous polynomial of degree $n-1$ and $\d_{\sD}^{ls}$ the partial holomorphic connection on $N_{\sD}^M$ defined as in Proposition \ref{LS_conn_new}. Let $\n_{\sD}^{ls}$ be any $(1,0)$-type extension of it and  $\n_1^0$ and $\n_2^0$ some connections respectively on $N_{S^0}^{\ot\nu}$ and $TM|_{S^0}$ such that the triple $(\n_1^0,\n_2^0,\n_{\sD}^{ls})$ is compatible with the short exact sequence
\begin{equation}\label{LS_seq}
0\longrightarrow N_{S^0}^{\ot\nu}\stackrel{\sD|_{S^0}}{\longrightarrow}TM|_{S^0}\stackrel{\r}{\longrightarrow} N_{\sD}^M\longrightarrow 0.
\end{equation}
Then set $\n_{\sD}^{\bullet}=(\n_1^0,\n_2^0)$ and observe that 
$$
\f\left(\n_{\sD}^{\bullet}\right)=\f\left(\n_{\sD}^{ls}\right)=0
$$
by the compatibility of the triple and the Bott vanishing theorem. This is the reason behind these choices. Let now $\S$ be a connected component of $\text{Sing}(S)\cup\text{Sing}(\sD)$ and $U\subset M$ an open subset such that $U\cap (\text{Sing}(S)\cup\text{Sing}(\sD))=\S$. Choose any connections $\n_1^U$ and $\n_2^U$ on $[S]|_U^{\ot\nu}$ and $TM|_U$ and set $\n_U^{\bullet}=(\n_1^U,\n_2^U)$. Finally,  let $\tilde{R}\subset U$ be any compact real sub-manifold of dimension $2n$ oriented as $M$ and such that $\S\subset\text{int}(\tilde{R})$. Assume its boundary $\partial\tilde{R}$ transverse to $S$ and with the orientation induced by $\tilde{R}$. Set $R=\tilde{R}\cap S$ and $\partial R=\partial\tilde{R}\cap S$. Then by definition the residue is 
$$
\text{Res}\left(\sD;\left.TM\right|_S-\left.[S]\right|_S^{\ot\nu};\S\right)=\int_R\f\left(\n_U^{\bullet}\right)-\int_{\partial R}\f\left(\n_{\sD}^{\bullet},\n_U^{\bullet}\right),
$$
where $\f(\n_{\sD}^{\bullet},\n_U^{\bullet})$ is the Bott difference form, defined as at \cite[pp.71-72]{Su1998}. One can show that this formula does not depend on the choices done.\par
In particular if $\S=\{p\}$ is an isolated point we can shrink $U$ so that $[S]^{\ot\nu}|_U$ and $TM|_U$ are trivial and we can assume $\n_1^U$ and $\n_2^U$ trivial respect to some local frames. Then $\f(\n_U^{\bullet})= 0$ and the formula reduces to
$$
\text{Res}\left(\sD;\left.TM\right|_S-\left.[S]\right|_S^{\ot\nu};p\right)=-\int_{\partial R}\f\left(\n_{\sD}^{\bullet},\n_U^{\bullet}\right).
$$
Assume now $p\in\text{Sing}(\sD)-\text{Sing}(S)$ and let $(U,z)\in\gU$ be a local chart at $p$. Observe that $N_{S^0}^{\ot\nu}$ and $TM|_{S^0}$ are not in general $N_{S^0}^{\ot\nu}$-bundles but 
$$
\left.N_{S^0}^{\ot\nu}\right|_{U\cap S^0}\qquad\text{ and }\qquad \left.TM\right|_{U\cap S^0}
$$
are canonically $N_{U\cap S^0}^{\ot\nu}$-bundles, thanks to the natural `holomorphic actions'  of $X$ on them given respectively by the Lie brackets 
$$
[X,\cdot]\qquad\text{ and }\qquad\left.[\cX,\cdot]\right|_{U\cap S^0}.
$$
These actions induce partial holomorphic connections on $N_{S^0}^{\ot\nu}|_{U\cap S^0}$ and $TM|_{U\cap S^0}$ along $N_{U\cap S^0}^{\ot\nu}$ which are, together with $\d_{\sD}^{ls}$, compatible with (\ref{LS_seq}) restricted to $U\cap S^0$. Therefore, similarly to Section \ref{BaumBott}, we can assume $\n_1^0$ and $\n_2^0$ to be `$X$-connections' defined on $U\cap S^0$ (and such that the triple $(\n_1^0,\n_2^0,\n_{\sD}^{ls})$ is compatible with (\ref{LS_seq})), and one can show that in this case
$$
\text{Res}\left(\sD;\left.TM\right|_S-\left.[S]\right|_S^{\ot\nu};p\right)=-\int_{\partial R}\f\left(\n_2^0,\n_2^U\right).
$$
Thus we can work again as in the proof of \cite[Th.III.5.5.]{Su1998} and obtain a formula like the one in \cite[Th.IV.5.3.]{Su1998}. In particular, letting the $h^j\in\cO_{M,p}$ be as in (\ref{germs_h}) and taking $\n_2^U$ trivial respect to the local frame  $\{\frac{\partial}{\partial z^1},\dots,\frac{\partial}{\partial z^n}\}$ of $TM$, we get the formula
\begin{equation}\label{LS_formula_1}
\text{Res}\left(\sD;\left.TM\right|_S-\left.[S]\right|_S^{\ot\nu};p\right)=\int_{\G}\frac{\f\left(-H\right)}{h^2\cdots h^n}\text{ d}z^2\wedge\cdots\wedge\text{d}z^n,
\end{equation}
where 
$$
H=\left(\left. \frac{\partial h^j}{\partial z^k}\right|_{U\cap S}\right)_{j,k=1,\dots, n}
$$
and $\G$ is the $(n-1)$ cycle defined as in Section \ref{BaumBott}. Observe that (\ref{LS_formula_1}) is very similar to (\ref{BB_formula_1}), there is just a little difference between the two matrices $H$.\par
\begin{remark}
Implicitly we have proved (and used) that
$$
\text{Res}\left(\sD;\left.TM\right|_S-\left.[S]\right|_S^{\ot\nu};p\right)=\text{Res}_{\f}\left(\cX;TM|_S;p\right) 
$$
where the residue on the right is the one associated to the natural `holomorphic action' of $X$ on $TM|_{U\cap S^0}$ mentioned just above (see also \cite[Rmk.IV.5.7.(1)]{Su1998} or \cite[Rmk.6.3.2(1)]{BSS2009}). In particular, we could have used directly the formula of \cite[Th.IV.5.3.]{Su1998}.  
\end{remark}
Lastly, assume $p\in\text{Sing}(S)$. In this case we can work exactly as in the last part of Section \ref{CamachoSad}. Assume to have an open neighborhood $U\subset M$ of $p$ such that $U\cap(\text{Sing}(S)\cup\text{Sing}(\sD))=\{p\}$ and $g|_U$ is a biholomorphism onto its image. Suppose moreover to have a local generator $y$ of $\cI_S$ and coordinates $(u^1,\dots,u^n)$ on it, with $(v^1=u^1\circ g^{-1},\dots,v^n=u^n\circ g^{-1})$ the corresponding ``special'' coordinates on $g(U)$ at $f(p)$. Let $\cV_{f,g}$ be the holomorphic vector field on $U$ defined as in (\ref{vector_V}), then reminding (\ref{LS_conn_loc}) and that $\nu>1$, equation (\ref{vectors_V_X}) implies also that the natural `holomorphic action' of $\cV_{f,g}$ on $N_{\sD}^M|_{U\cap S^0}$ by Lie bracket induces locally $\d^{ls}_{\sD}$.  Hence if we choose the coordinates $(u^1,\dots,u^n)$ in such a way that $\{y=(v^2\circ f-v^2\circ g)/y^{\nu}=\cdots=(v^n\circ f-v^n\circ g)/y^{\nu}=0\}=\{p\}$  then we can express the residue at $p$ by the formula in \cite[Th.IV.5.3.]{Su1998}. In particular, we get the formula
\begin{equation}\label{LS_formula_2}
\text{Res}\left(\sD;S;p\right)= \int_{\G} \frac{y^{\nu(n-1)}\f(-Y)}{\prod_{j=2}^n(v^j\circ f-v^j\circ g)} \text{ d}u^2\wedge\cdots\wedge\text{d}u^n,
\end{equation}
where $Y$ is the Jacobian matrix of 
$$
\left((v^2\circ f-v^2\circ g)/y^{\nu},\dots,(v^n\circ f-v^n\circ g)/y^{\nu}\right)
$$
with respect to the coordinates $(u^1,\dots,u^n)$ and $\G$ is the $(n-1)$ cycle
$$
\G=\left\{q\in U\cap S\text{ s.t. }\left|\frac{v^j\circ f-v^j\circ g}{y^{\nu}}(q)\right|=\e, \text{ for }j=2,\dots,n\right\},
$$
for $\e>0$ small enough, oriented as usual.\par
\bigskip

 \end{document}